\documentclass{article}
\usepackage[utf8]{inputenc}
\usepackage{lineno,hyperref}
\usepackage{amsmath,amssymb}
\usepackage{graphicx}
\usepackage{graphics}
\usepackage[english]{babel}
\usepackage{comment} 
\modulolinenumbers[10]
\modulolinenumbers[1]
\usepackage{amsthm}
\usepackage{color}
\usepackage{geometry}		
\usepackage{mathtools}
\usepackage{caption}
\usepackage{subcaption}
\usepackage{lineno}
\usepackage{multirow}

\modulolinenumbers[10]
\modulolinenumbers[1]
\newcommand\nnfootnote[1]{%
  \begin{NoHyper}
  \renewcommand\thefootnote{}\footnote{#1}%
  \addtocounter{footnote}{-1}%
  \end{NoHyper}
}

\newcommand*{\blue}{\textcolor{black}}
\newcommand*{\blu}{\textcolor{black}}
\title{Dynamical effects of electromagnetic flux on Chialvo neuron map: nodal and network behaviors}
\begin{document}
\author{Sishu Shankar Muni \\
Department of Physical Sciences,\\
Indian Institute of Science and Educational Research Kolkata,\\
Campus Road, Mohanpur, West Bengal, 741246,\\
India\\
\newline
School of Mathematical and Computational Sciences,\\
Massey University,\\
Colombo road, Palmerston North, 4410,\\
New Zealand\\
\newline\\
Hammed Olawale Fatoyinbo, Indranil Ghosh\\
School of Mathematical and Computational Sciences,\\
Massey University,\\
Colombo road, Palmerston North, 4410,\\
New Zealand\\}
\maketitle
\begin{abstract}
We consider the dynamical effects of electromagnetic flux on the discrete Chialvo neuron. It is shown that the model can exhibit rich dynamical behaviors such as multistability, firing patterns, antimonotonicity, closed invariant curves, various routes to chaos, fingered chaotic attractors. The system enters chaos via period-doubling cascades, reverse period-doubling route, antimonotonicity, via closed invariant curve to chaos. The results were confirmed using the techniques of bifurcation diagrams, Lyapunov exponent diagram, phase portraits, basins of attraction and numerical continuation of bifurcations. Different global bifurcations are also shown to exist via numerical continuation. After understanding a single neuron model, a network of Chialvo neuron is explored. A ring-star network of Chialvo neuron is considered and different dynamical regimes such as synchronous, asynchronous, chimera states are revealed. Different continuous and piecewise continuous wavy patterns were also found during the simulations for negative coupling strengths.  
\end{abstract}

\nnfootnote{\blue{Keywords: Chialvo neuron map, multistability, routes to chaos, antimonotonicity, fingered attractors, closed invariant curves, ring-star network, chimera states, traveling waves.}}
\nnfootnote{E-mail addresses: S.S. Muni (\url{sishu1729@iiserkol.ac.in, ssmuni760010@gmail.com}), H. O. Fatoyinbo (\url{h.fatoyinbo@massey.ac.nz}), I. Ghosh(\url{i.ghosh@massey.ac.nz})}

\section{Introduction}
Neurons are important dynamical units of the brain. Since neurons fire potentials and behave dynamically, it is considered a dynamical system. For many years, researchers have considered many dynamical system models to understand the behavior of neurons. Hodgkin-Huxley model was the first neuron model developed and a lot of firing patterns, bursting patterns were shown which mimic the real-time behavior of neurons \cite{HoHu52}. Since then many continuous and discrete neuron dynamical systems have been formulated and explored \cite{FHN61, UsSu19, ML81}. Recently, effects of external electromagnetic field on neurons have been considered \cite{em17, Muni22a}. The motivation behind this is that under the action of electromagnetic field on embryonic neurons has led to their growth in the cultures. In \cite{HGu13}, researchers have shown a very good agreement between the dynamical behaviors of a theoretical Hindmarsh-Rose neuron model and experimental neuron firing patterns. Experimental studies have been done on firing patterns, and bifurcations that resembled those in the physiologically based Chay neuron model \cite{GuPaCh14}.  A new type of experiment has been done which reveals that a neuron functions as multiple threshold units in \cite{Sar17}. Basic bifurcation structures are eminent in injured nerve fibers in \cite{BJia17}. They showed a basic bifurcation structure from bursting to spiking in two-dimensional parameter space. Experimental impulse recording from hypothalamic brain slicing has revealed an important transition from tonic firing to bursting discharges studied in \cite{Braun11}. Chialvo neuron map is a discrete neuron model which is introduced by Chialvo \cite{Ch95}. There have been researches showing quasiperiodic structures in the map, analytically studying the Chialvo map \cite{Jh06,Wa18}. Chialvo neuron map accounts for the dynamics in perturbed excitable systems, phase-locked responses in parameter space. This shows that the Chialvo map is generic to a large class of neuron systems. Especially discrete neuron maps are of interest to understanding complex dynamics arising in large distributed neuron media. This paper improves the Chialvo neuron model by studying the effects of electromagnetic flux. 

\blue{Various complex dynamical behaviors of a single neuron or a network of neurons are important to understand the complex behaviors in the \blue{brain like serious neurological diseases} \cite{IZHBook}. One of the main attributes of neuron systems is information processing. It depends not only on the electrophysiological properties but also on the dynamical properties of the neurons. Therefore it is of utmost importance to uncover the dynamical behaviors, bifurcation patterns, dynamical states of neuron cells. Bifurcation theory determines the excitable properties of the neuron. This motivates to explore the complex dynamical behaviors of the neuron systems.}

The aim of this research is \blu{threefold}. First is to consider the effect of external electromagnetic field on the Chialvo neuron and explore its dynamical behavior, second is to consider the network of Chialvo neurons under electromagnetic field and observe their behavior, and \blue{finally} numerically continuating bifurcations in order to detect global bifurcations. This is more natural as after \blu{analyzing} a single neuron, we proceed to explore a network of neurons as in the \blu{natural} nervous system, the neurons are  interconnected in a very complex \blue{topological fashion}. We have shown that the dynamical behavior of Chialvo \blu{neurons} under the action of \blu{the} electromagnetic field is itself rich. Multistability in Chialvo map under the action of electromagnetic flux, various global bifurcation \blu{phenomena}, three different routes to chaos (\blu{period-doubling} route to chaos, antimonotonicity route to chaos via periodic and chaotic bubbles, reverse bubbles, and via an invariant closed curve to chaos), chaotic fingered attractors, chimera states in a ring-star network of Chialvo neurons are firstly studied in this research \blue{which build-up to} the novelty of this paper. 

Multistability is an exotic phenomenon in the theory of dynamical systems. It refers to the coexistence of several attractors at the same set of parameter values. Multistability has been found in many \blu{real-world} systems like visual optical illusions, convection currents, electronic circuits, \blu{and} engineering systems \cite{Feu18}. It is undesirable in many engineering systems as \blu{the} presence of several states might lead to the malfunctioning of systems. \blu{In contrast}, it is an advantageous feature for the information processing of neurons. Bistability is a common feature observed in many neuron models like the Hodgkin-Huxley model \cite{Kam18}, cerebellar Purkinje cells \cite{Jo13} \textcolor{black}{, and muscle cells \cite{hammed}}. Multistability also \blu{favors} the multifunctionality of central pattern generators \cite{Muni22a}. \blu{The} coexistence of periodic attractors of different periods is shown in the case of \blu{the} Chialvo neuron \blu{model} in this paper. To aid the multistable phenomenon, we have illustrated it through the basins of attraction diagram.  

\textcolor{black}{Bifurcation analysis is another tool used to investigate the variety of possible dynamical \blu{behaviors} in neurons. Transitions between different states, such as from quiescent to periodic, and from periodic to complex periodic in neuron models, have been widely studied \cite{Tsumoto2006BifurcationsModel, Duan2006Codimension-twoModel,WU2021,hammed1}. Jing {\em et al}~\cite{Jh06}, studied the influence of different parameters on the dynamical \blu{behavior} of the Chialvo map via bifurcation analysis. They showed several bifurcation structures and codimension-1 bifurcations that \blu{characterized various behaviors,} such as period-doubling and inverse period-doubling bifurcations, in the Chialvo map. Also, Wang and Cao \cite{Wa18} confirmed the existence of mode-locking and quasiperiodicity behaviour in some parameter regions of the Chialvo map. However, in this work, we consider one- and two-parameter bifurcation analysis to study the effects of electromagnetic flux on the dynamics of the Chialvo map. We show some global bifurcations that have not been reported in previous studies \blu{of} the Chialvo map.}

Since Huygen's classical work on coupled pendulums \cite{Huy17}, many oscillators were coupled, and researchers have \blu{analyzed} and studied their behavior. They found that oscillators synchronize and desynchronize under certain conditions. Later Kuramoto discovered that a third kind of dynamical state is also possible in the networks of coupled oscillators known as chimera states. Various \blu{network elements} have similar \blu{phases and frequencies} giving rise to the phenomenon of synchronization. Sometimes, a network with synchrony can pass through an intermediate state called a chimera, where there is a coexistence of synchronous and asynchronous states \cite{kuramoto2002coexistence, panaggio2015chimera}.  
\blu{Many} researchers have focused on studying chimera states in different dynamical systems \blu{to} date \cite{scholl2016synchronization, majhi2019chimera, omel2018mathematics}. Chimera states are also found in \blu{real-world} scenarios such as dolphin's sleeping pattern \cite{glazeneural} with one eye open and closed (a part of \blu{the} brain is active while the other is at rest), flashing of fireflies \cite{haugland2015self}, social systems and many more \cite{buscarino2015chimera, martens2010solvable, zakharova2016amplitude, omelchenko2015robustness, hizanidis2014chimera}. Scientists have even uncovered epilepsy \blu{and} \blu{schizophrenia} as topological diseases \cite{uhlhaas2006neural} that is it depends on the topology of the neurons connected. This motivates us to study chimera states in a network of Chialvo neurons. 

In this paper, we have considered \blu{the} ring-star network introduced in \cite{Muni20} of Chialvo neurons and have explored the collective behavior of the neurons with variation in the coupling parameters. Efforts have been made to understand the \blu{neuron's behavior} under different network topologies formed by considering different combinations of coupling parameters. Star and ring network topologies have \blu{typical} applications in social systems, social networks, computer networks \cite{roberts1970computer}, etc.  \blue{The ring-star topology of neurons has been considered to model our system because we believe that it acts as a trade-off between both the complexity and \blue{the reality of a nervous system}. The ring-star topology, in some sense, imitates the interconnections among the neurons through chemical and electrical synapses, which \blu{are} essentially responsible for transfer of information among them. We also notice that the topology does not fail to capture \blu{actual world} occurrences like synchronization and chimera.}

The paper is \blue{organized} as follows. \S~\ref{sec:actionEM} introduces the model and improves it under the action of electromagnetic flux. \S~\ref{sec:features} discusses the features of the improved Chialvo map, such as fixed points \blu{and} non-invertibility. \S~\ref{sec:multistable} shows the exotic phenomenon of multistability. \S~\ref{sec:bifstruc} illustrates the bifurcation analysis in Chialvo neurons. \S~\ref{sec:numbif} discusses continuation of bifurcations, two parameter bifurcation analysis using software {\sc MatContM}. \S~\ref{sec:firing} discusses the bursting and spiking patterns \blu{exhibited} by the Chialvo neuron map under the action of electromagnetic flux with the variation of parameters. \S~\ref{sec:chimera} introduces a ring-star network of Chialvo neurons and explores the spatiotemporal patterns, chimera states in the network under \blue{both positive and negative coupling strengths}. \S~\ref{sec:conclusions} provides some future directions of this research and concludes the paper. All the simulations in the paper were carried out using {\sc MATLAB} and Python.


\section{Chialvo map under the action of electromagnetic flux coupling}
\label{sec:actionEM}
\begin{figure}[!htbp]
\begin{center}
\includegraphics[width=0.7\textwidth]{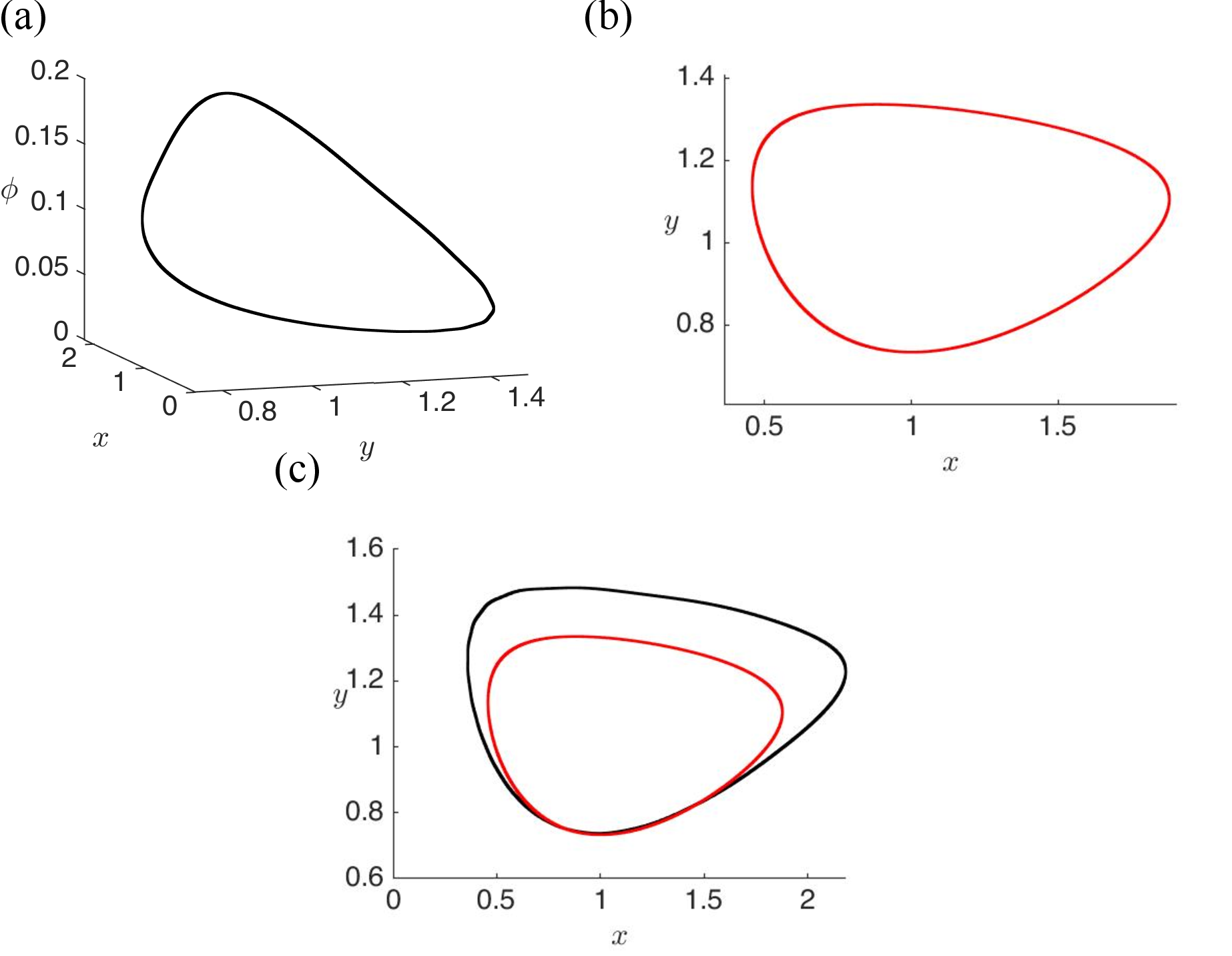}
\end{center}
\caption{\blue{ In (a), a three dimensional closed invariant curve in black is shown in the presence of electromagnetic flux $\phi$. In (b), a two-dimensional closed invariant curve is shown in red in the absence of electromagnetic flux. In (c), an overlap plot of the attractors in the absence (red) and presence (black) of electromagnetic flux is considered. It is found that the attractors are indeed distinct. The parameters considered are $a=0.89, b=0.18, c=0.28, k_{0}=0.06$. With the inclusion of electromagnetic flux parameters $k=-0.2, \alpha = 0.1, \beta = 0.2, k_{1}=0.1, k_{2}=0.2$. }} 
\label{fig:DynCompare}
\end{figure}

Chialvo introduced a two dimensional discrete neuron system given by

\begin{equation}
\begin{aligned}
    x_{n+1} &= x_{n}^{2} e^{(y_{n} - x_{n})} + k_{0},\\
    y_{n+1} &= a y_{n} - b x_{n} + c
\end{aligned}
  \label{eq:ChialvoMap}
\end{equation}
where $x$ denotes the activation or potential variable, \blue{and} $y$ denotes the recovery like variable. The parameters are $a,b,c,k_{0}$ where $a$ represents the time constant of recovery $(a<1)$, $b$ denotes the activation dependence of the recovery process $(b<1)$, $c$ denotes the offset, and $k_{0}$ the parameter which can act either as a constant bias or as a time dependent additive perturbation. In this paper, we have considered $k_{0}$ to be time independent parameter. The Chialvo neuron is now improved considering the effect of electromagnetic flux. A memristor is usually effective in describing the effects of electromagnetic flux in neuron systems.  From a physical point of view, it bridges the gap between the effect of describing the electromagnetic induction and membrane potential of the neuron.  Due to the electromagnetic flux, we have an induced induction current given by 
\begin{equation}
    \frac{dq(\phi)}{dt} = \frac{dq(\phi)}{d\phi} \frac{d\phi}{dt} = M(\phi) \frac{d\phi}{dt} = kM(\phi)x
    \label{eq:indcurr}
\end{equation}
We see that the effect of electromagnetic induction and electric field can be described via induction current. Here $\phi$ denotes the magnetic flux across the neuron membrane, $k$ represents the strength of electromagnetic flux coupling and $M(\phi)$ denotes the memconductance of magnetic flux controlled memristor. \textcolor{black}{We note that $k$ can take both positive and negative values.} Researchers have used many forms of memconductance function and it was mentioned in \cite{Ch12} about the importance of memristors in understanding about the action potential in neurons. In this study we take the commonly used memdconductance function given by $M(\phi) = \alpha + 3 \beta \phi^2$. Thus
under the action of electromagnetic flux, the Chialvo neuron in \eqref{eq:ChialvoMap} can be improved to the following map 
\begin{equation}
\begin{aligned}
\label{eq:ChialvoMag}
x_{n+1} &= x_{n}^{2} e^{(y_{n}-x_{n})} + k_{0} + k x_{n} M(\phi_{n}),\\
y_{n+1} &= a y_{n} - bx_{n} + c,\\
\phi_{n+1} &= k_{1}x_{n} - k_{2}\phi_{n},
\end{aligned}
\end{equation}
where an additional induction current term due to electromagnetic flux is added in the $x$ state variable and this makes the  Chialvo map to be a three dimensional smooth map. \textcolor{black}{In \eqref{eq:ChialvoMag}, the $k_{1}x$ term denote the membrane potential induced changes on magnetic flux and the $k_{2}\phi$ term denote the leakage of magnetic flux \cite{WuLevy17}}. The electromagnetic flux parameters constitute $\alpha, \beta, k_{1}, k_{2}$. Observe that the map is asymmetric in the sense that the form of the map changes under the transformation $(x,y,\phi) \rightarrow (-x,-y,-\phi)$. From a dynamical system theory perspective, the improved three dimensional Chialvo map in \eqref{eq:ChialvoMag} is an important system as it can show many new kinds of global bifurcations which cannot be exhibited by two dimensional map.  An important aspect of this paper is to study the effect of electromagnetic flux on the neurons and we observe the behavior of the system with the variation of the electromagnetic coupling strength $k$.
\blue{Since electromagnetic flux has physiological effects on the neuron system and a time varying flux is considered, the electromagnetic flux $\phi$ is considered as a  state variable and further improves the Chialvo neuron model transforming it to a three-dimensional map.  To observe the effect of the flux $\phi$ on the Chialvo neuron transforming it to a three dimensional map, a comparision between the two dimensional and improved three dimensional map is shown in Fig. \ref{fig:DynCompare}. In Fig. \ref{fig:DynCompare}(b), we observe a closed invariant curve in red and when flux $\phi$ is introduced, the three dimensional neuron map exhibits a three dimensional closed invariant curve in black as shown in Fig. \ref{fig:DynCompare}(a). In Fig. \ref{fig:DynCompare}(c), the projection plot is considered for both the three dimensional and two dimensional closed invariant curve and it is found that they are distinct and don't overlap.  }
In the next section, we analyse the fixed points of the improved Chialvo map in \eqref{eq:ChialvoMag}.
\section{Features of Chialvo neuron}
\label{sec:features}
\subsection{Fixed points of the proposed system}
\begin{figure}[!htbp]
\begin{center}
\includegraphics[width=1\textwidth]{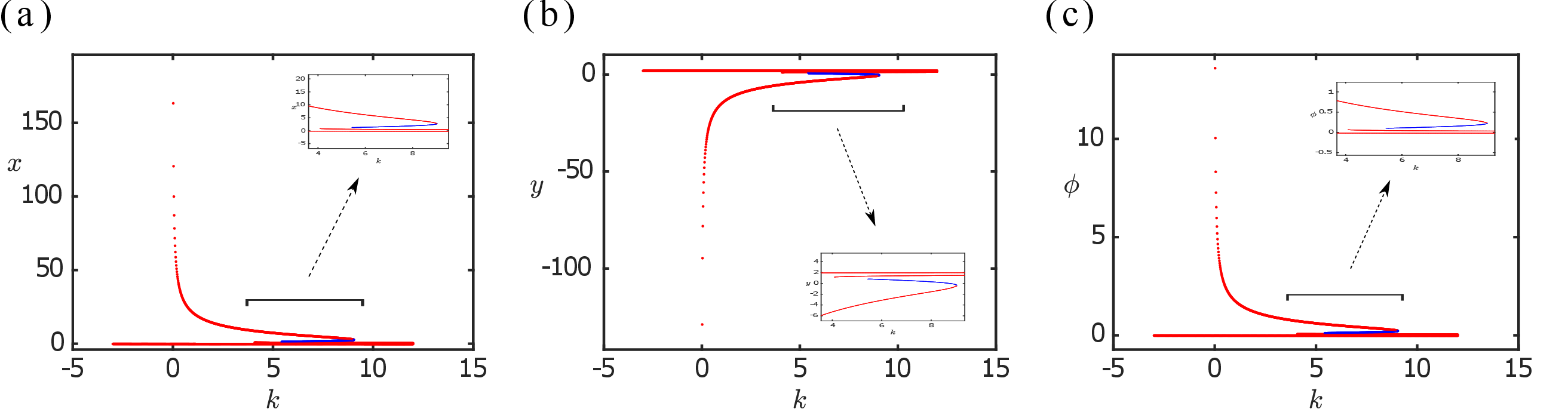}
\end{center}
\caption{\blue{ Bifurcation diagram of fixed points with respect to the parameter $k$. In (a), $x$ component of the fixed point is \blue{shown}. In (b), $y$ component of the fixed point is shown. In (c), $\phi$ component of the fixed point is shown. A fixed point is marked by red colour if it is a saddle and is marked by blue if it is an asymptotically stable fixed point. Parameter $k$ is varied in $[-3,12]$. The remaining parameters are set as $a=0.5, b=0.4, c=0.89, k_{0}=-0.44, k_{1}=0.1, k_{2}=0.2, \alpha = 0.1, \beta = 0.1.$}} 
\label{fig:fpKChialvo}
\end{figure}

\begin{figure}[!htbp]
\begin{center}
\includegraphics[width=0.5\textwidth]{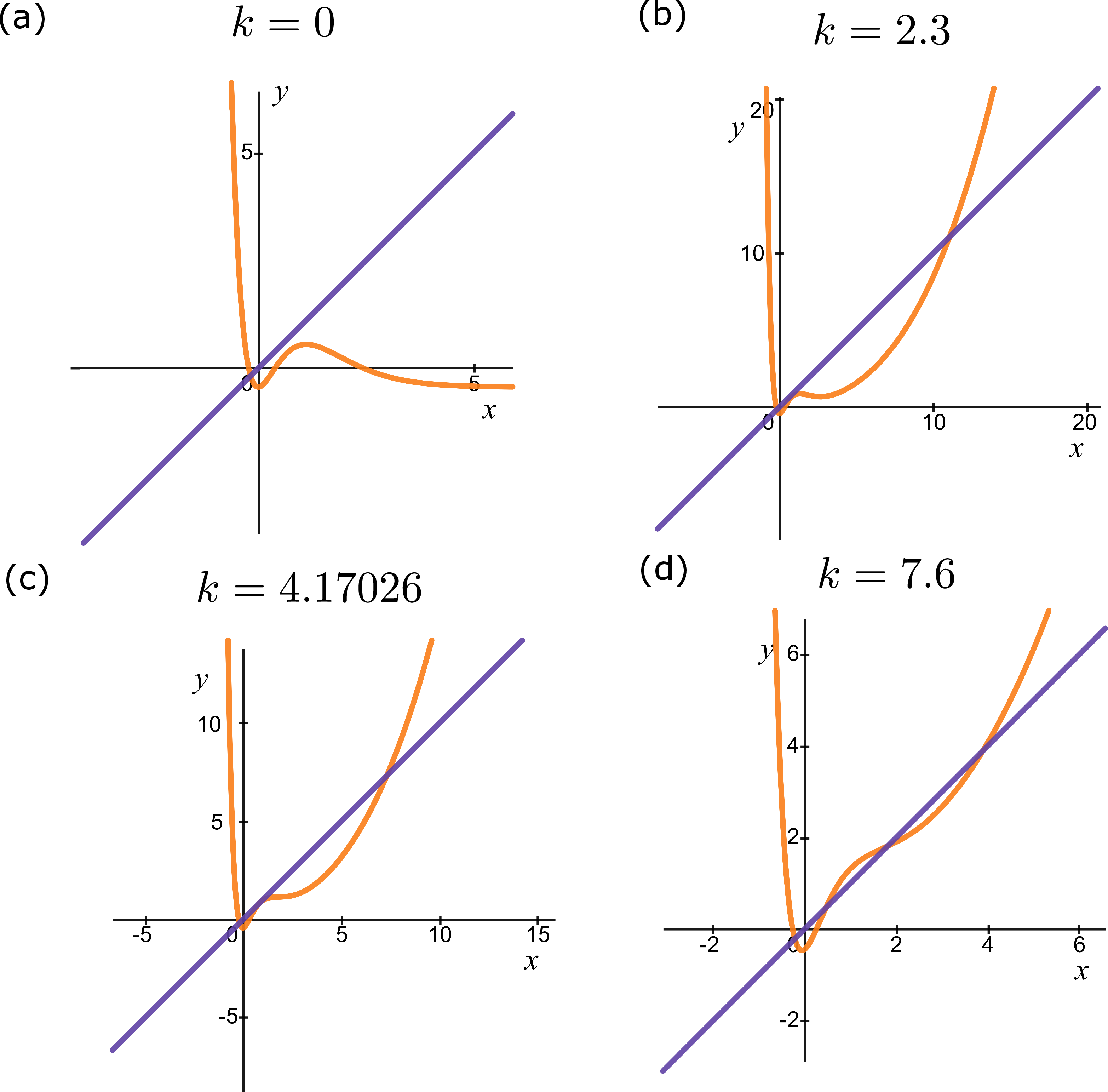}
\end{center}
\caption{\blue{Number of fixed points varies as the electromagnetic flux parameter $k$ is varied. The parameters are set as $a=0.5, b=0.4, c=0.89, k_{0}=-0.44, k_{1}=0.1, k_{2}=0.2, \alpha = 0.1, \beta = 0.1.$}} 
\label{fig:fpChialvo}
\end{figure}

A fixed point analysis of the two-dimensional Chialvo map has been done in \cite{Jh06}. In this section, we attempt to understand the number of fixed points exhibited by system \eqref{eq:ChialvoMag}. Researchers have made progress \blu{in analyzing} the fixed points in the case of \blu{the} linear three-dimensional \blu{maps} in \cite{Ra17}. The map in \eqref{eq:ChialvoMag} involves exponential terms, and hence the analytical study of the fixed point gets difficult. The fixed point \blue{of} \eqref{eq:ChialvoMag} \blue{are} given by solving the three simultaneous equations
\begin{equation}
    \label{eq:SimulFP}
    \begin{aligned}
    x^2 e^{y-x} + k_{0} + kx(\alpha + 3 \beta \phi^2) &= x,\\
    ay-bx+c &= y,\\
    k_{1}x - k_{2}\phi &= \phi.
    \end{aligned}
\end{equation}
This leads to solving a transcendental equation 
\begin{equation}
x^{2} e^{\frac{(b-a+1)x -c}{a-1}} + k_{0} + \frac{3k \beta k_{1}^{2}}{\blue{(1+k_{2})^{2}}}x^3 + x k \alpha = x
\label{eq:transcendental}
\end{equation}
Equation \eqref{eq:transcendental} can be solved numerically or graphically. We have used graphical method to solve the transcendental equation. It can be thought of as intersection of two curves $y=f(x)$ and $y=g(x)$, where $f(x) = x^{2} e^{\frac{(b-a+1)x -c}{a-1}} + k_{0} + \frac{3k\beta k_{1}^{2}}{\blue{(1+k_{2})^{2}}}x^3 + x k\alpha$ and $g(x) = x$. \blue{To explore the variation of the fixed points with respect to the electromagnetic flux parameter $k$, a bifurcation diagram of fixed points with respect to parameter $k$ is shown in Fig. \ref{fig:fpKChialvo}. A red point denote a saddle fixed point and a blue point denote an asymptotically stable fixed point. The number of fixed points increases with increase in parameter $k$. Observe that for $5<k<10$, the saddle and asymptotically stable fixed point come closer and collide. This diagram also sheds light on the number of fixed points with the variation of $k$. With few discrete values of $k$, the graph of fixed points \blue{are} shown in Fig. \ref{fig:fpChialvo}.}

This is shown in Fig.~\ref{fig:fpChialvo}, where the orange curve denotes $y=f(x)$ and the straight line in violet denotes $y=x$.
The fixed points of \eqref{eq:ChialvoMag} are given by $(\tilde{x},\tilde{y},\tilde{\phi})$ where $\tilde{x}$ satisfies \eqref{eq:etranscendental},  $\tilde{y}$ is given by \eqref{eq:ycompFP}, $\tilde{\phi}$ is given by \eqref{eq:phicompFP}.
\begin{equation}
\tilde{x}^{2} e^{\frac{(b-a+1)\tilde{x} -c}{a-1}} + k_{0} + \frac{3k\beta k_{1}^{2}}{\blue{(1+k_{2})^{2}}}\tilde{x}^3 + \tilde{x} k\alpha = \tilde{x}
\label{eq:etranscendental}
\end{equation}

\begin{equation}
    \tilde{y} = \frac{b\tilde{x}-c}{a-1}
    \label{eq:ycompFP}
\end{equation}

\begin{equation}
    \tilde{\phi} = \frac{k_{1}x}{1+k_{2}}
    \label{eq:phicompFP}
\end{equation}

The number of fixed points depends on the number of real solutions to \eqref{eq:etranscendental}. Observe that number of fixed points also depends on the parameters. Below we show that the system possesses many fixed points as parameters are varied and showcase the stability and type of fixed point. For simplicity, we vary a single parameter in the system and observe how the number and \blu{the} stability of fixed point change. Here we vary the electromagnetic flux parameter $k$ and showcase that the system can exhibit one, two, three, and four fixed points. We also \blu{analyze} the stability of the fixed points via the eigenvalues. Such an analysis is presented in table \ref{Tab:stabFixed} in accordance with Fig. \ref{fig:fpChialvo}.

\begin{table}[htbp]
\centering
\caption{\blue{Stability analysis of the fixed points}}
\label{Tab:stabFixed} 
\begin{tabular}{ |c|c|c|c|  }
\hline
$k$ & ($x, y, \phi$) & Eigenvalues ($\lambda_1, \lambda_2, \lambda_3$) & Type \\
\hline\hline
\hline
\multirow{1}{*}{0}&(-0.1787, 1.9230, -0.0149) & (-0.2,0.4714,-3.1566) & saddle\\
\hline
\hline
\multirow{2}{*}{2.3}&(-0.1883, 1.9306, -0.0157) & (-3.1669,  0.4678, -0.199855) & saddle\\
&(12.953, -8.5824, 1.0794) & (1.93686, -1.1029,  0.5) & saddle\\
\hline
\hline
\multirow{3}{*}{4.1703}&(-0.1964, 1.9371, -0.0164) & (-3.1903,  0.4647, -0.1997) & saddle\\
&(8.546, -5.0568, 0.7122) & (1.8094, -0.9579,  0.5) & saddle\\
&(0.831, 1.1152, 0.0693) & (1.1991,  1.0206 , -0.2059) & saddle\\
\hline
\hline
\multirow{4}{*}{7.6}&(-0.212, 1.9496, -0.0177) & (-3.2712,  0.4586, -0.1994) & saddle\\
&(0.461, 1.4112, 0.0384) & (2.4908,  0.61003, -0.2026) & saddle\\
&(1.755, 0.3760, 0.1462) & (0.7453+0.4697i,  0.7453-0.4697i,
        -0.2735) & \blue{asym. stable}\\
&(4.559, -1.8762, 0.3299) & (-0.6224,  1.4026,  0.5092) & saddle\\
\hline
\end{tabular}
\end{table}

\subsection{Non-invertibility}
Non-invertibility of a discrete map is an \blu{essential} feature \blu{because} non-invertibility makes stretching and folding action in discrete maps. This can \blu{even be} seen in \blu{the} logistic map, a well-known one-dimensional non-invertible map. \blu{ Depending on the system's dimensionality, critical lines, curves, and surfaces separate different numbers of preimages in the phase space.} If the map is one-dimensional, we obtain a critical line \blu{that} separates intervals with a \blu{different} number of preimages. When the map is two-dimensional, critical curves separate regions into distinct preimages. In the case of three-dimensional maps, instead, we get surfaces in the three dimensional phase space separating regions with distinct preimages. \blue{Much on the classifications of non-invertible maps, its dynamical implications can be found in \cite{Mi96,CM96}. \blu{However, many clear visualizations of the critical surfaces remain an area to explore in the future for higher-dimensional maps.} }

In line with the notations from Mira's book in \cite{CM96}, there are various classifications of smooth non-invertible maps. We follow the same classification here. 

\begin{figure}[!htbp]
\begin{center}
\includegraphics[width=0.7\textwidth]{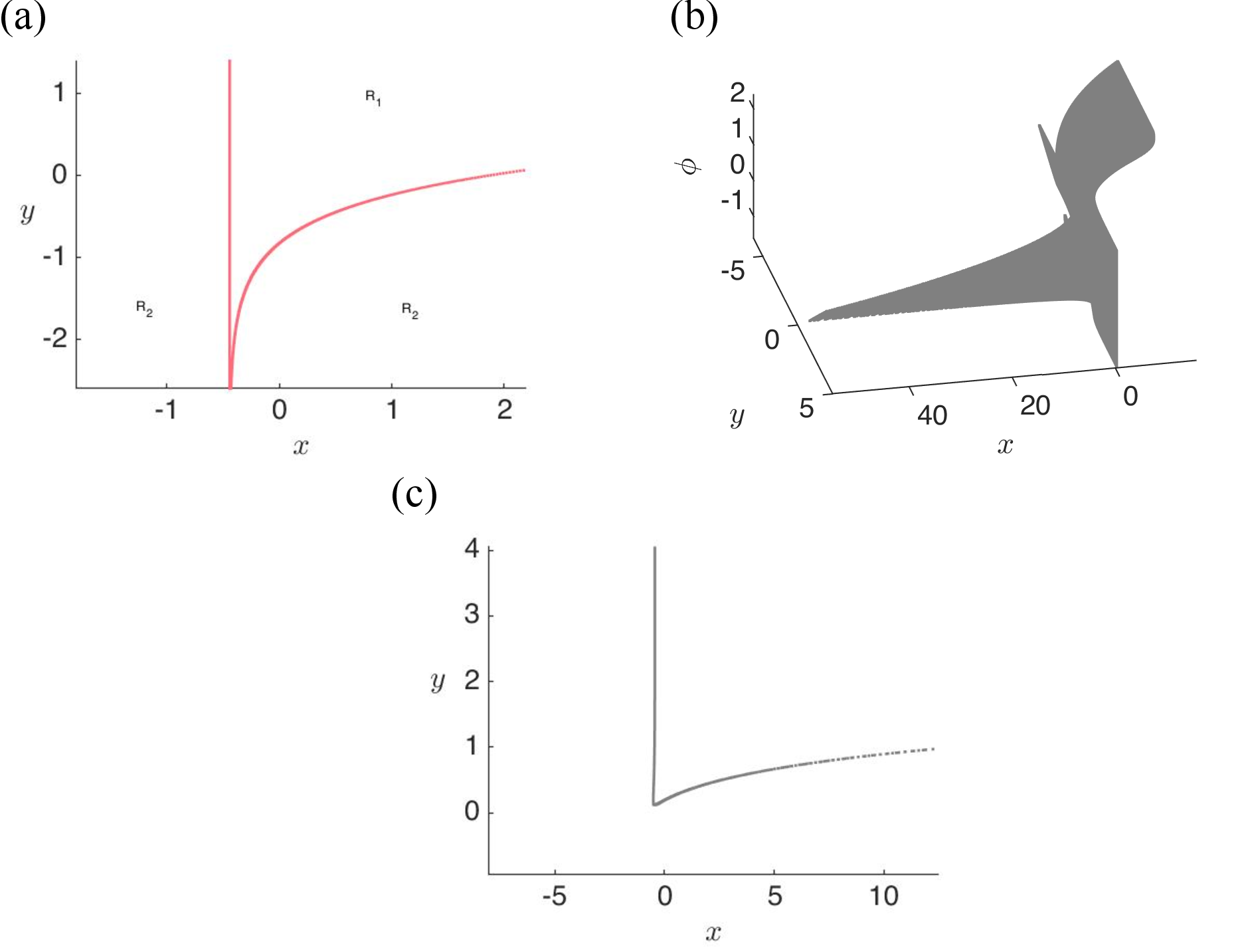}
\end{center}
\caption{In (a), the critical curve $LC$ is shown in red for the two-dimensional Chialvo map \eqref{eq:ChialvoMap}. Observe that the LC curve is a cusp separating the $xy$ region into two regions $R_{1}, R_{2}$. The region $R_{1}$ has three preimages whereas region $R_{2}$ has one preimage. In (b), the critical surface in $x-y-\phi$ space is \blue{shown} in \blue{grey for the three dimensional Chialvo map in \eqref{eq:ChialvoMag}}. \blue{A slice of the critical surface is shown in (c) for $\phi=3$, which is a cusp separating the phase space into two distinct regions with distinct preimages.} This also confirms that the three dimensional Chialvo map in \eqref{eq:ChialvoMag} is of type $Z_{1}-Z_{3}$. The parameters are set as $a = 0.5, b = 0.4, c = 0.89, k_{0} = -0.44, \alpha = 0.1, \beta = 0.1, k_{1} = 0.1, k_{2} =0.2, k=-0.1$.}
\label{fig:LCcurveNonInvertible}
\end{figure}

We show that the Chialvo map in \eqref{eq:ChialvoMap} and improved Chialvo map in \eqref{eq:ChialvoMag} are non-invertible. Let us analyse the non-invertibility of the two-dimensional Chialvo map. 
The Jacobian matrix of the Chialvo map \eqref{eq:ChialvoMap} is given by 
$$J_{1} = \begin{bmatrix} 2xe^{y-x} - x^2e^{y-x} & x^2 e^{y-x} \\
-b & a\\
\end{bmatrix}$$
The curve of merging rank-$1$ preimages is where the determinant of the Jacobian matrix vanishes and is given by  
\begin{equation}
    LC_{-1} = \{(x,y) | e^{y-x} \left(2ax - ax^2 + bx^2\right) =0\}
    \label{eq:LC-1Chialvo2D}
\end{equation}
The critical curve $LC$ (from the french word \textit{``Ligne Critique"}) is then given by $LC = f(LC_{-1})$. More generally critical curves divide the phase space into regions $R_{i}$ which has constant number of $k_{i}$ preimages. The map is then classified as $Z_{k_{1}}-Z_{k_{2}}-\ldots-Z_{k_{n}}$ as the type of regions $Z_{k_{i}}$ appear. It is found that \eqref{eq:ChialvoMap} is of the type $Z_{1}-Z_{3}$. It is shown after computing the critical curve $LC = f(LC_{-1})$ which is a cusp as shown in \blue{Fig. \ref{fig:LCcurveNonInvertible}}. This suggests that the cusp separates the region of three images on one side and one-preimage on the other side. \blue{We have also observed that the} subsequent images $f^{k}(LC_{-1})$ of the critical curve $LC_{-1}$, bound the attractor boundary. This possibly explains the shape of \blue{fingered attractors explored in the later sections of the paper}. 

It is next natural to think about the invertibility of the three-dimensional improved Chialvo map \eqref{eq:ChialvoMag}. The Jacobian matrix of the improved Chialvo map in \eqref{eq:ChialvoMag} is given by 
$$J = \begin{bmatrix} e^{y-x}(2x  -x^2) + k(\alpha + 3 \beta \phi^2) & x^2 e^{y-x} & 6kx\beta\phi\\
-b & a & 0\\
k_{1} & 0 & -k_{2}\end{bmatrix}.$$ The surface of merging rank-one preimages is given by 
\begin{equation}
     LC_{-1} = \{(x,y,\phi) |  -e^{y-x}(2x-x^2)k_{2} a - k_{2} a k(\alpha + 3\beta\phi^2) - bk_{2} x^2 e^{y-x} - 6kx\beta \phi a k_{1}=0\}.
    \label{eq:LC-1Chialvo3D}
\end{equation}
 The critical surface is then given by $LC = f(LC_{-1})$. The critical surface is shown in \blue{Fig.~\ref{fig:LCcurveNonInvertible} (b)}. \blue{We can observe a cusp in three dimensional surface. To visualise better, we can take a slice of the critical surface at some discrete value of $\phi$. Such a slice of the critical surface is shown in Fig.~\ref{fig:LCcurveNonInvertible}(c) for $\phi=3$. Observe that it is a cusp, separating the $x-y$ plane into two distinct regions of preimages similar to Fig.~\ref{fig:LCcurveNonInvertible}(a) confirming the three-dimensional Chialvo mp to be a non-invertible map of type $Z_{1}-Z_{3}$.}

\section{Multistability}
\label{sec:multistable}
\begin{figure}[!htbp]
\begin{center}
\includegraphics[width=0.9\textwidth]{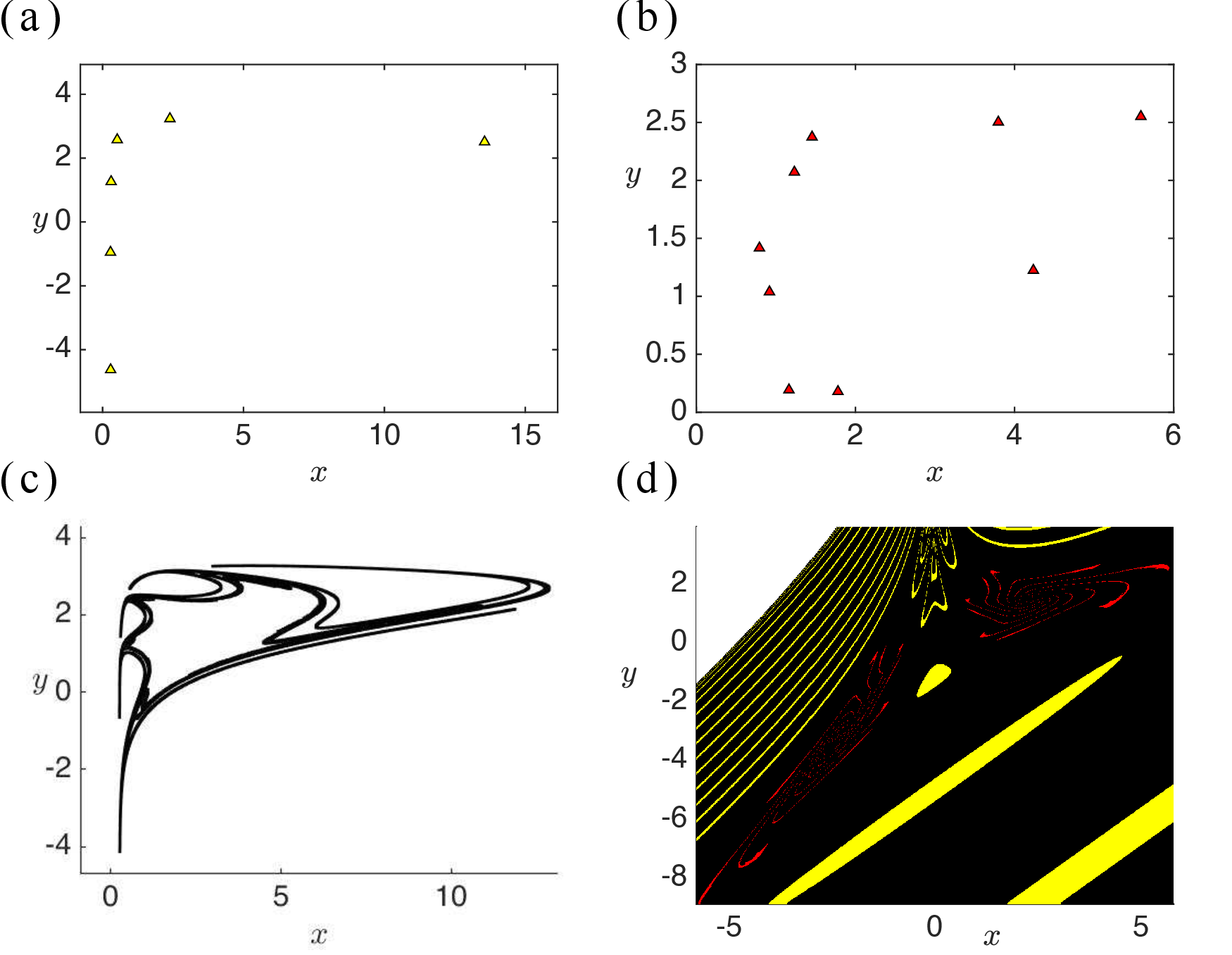}
\end{center}
\caption{\blue{Coexistence of a period-six, period-nine, chaotic attractor along with divergence. In (a), we observe a period-six orbit marked by yellow triangles. In (b), we observe a period-nine orbit marked by red triangles. In (c), we observe a chaotic attractor in black. In (d), basin of attraction diagram denoting the coexistence of stable period-six and period-nine solution with the chaotic attractor is shown. Yellow region corresponds to the region where the initial conditions converge to a period-six orbit, red colour where the initial conditions converge to a \blue{period-nine} orbit, black regions where the initial conditions converge to a chaotic attractor, white coloured regions where the trajectory diverges to infinity. }  \blue{The parameters are set as $a=0.6,b=0.6,c=2,k_{0}=0.28,k_{1}=0.1,k_{2}=0.2,k=0.002,\alpha=0.1, \beta=0.2$.}}
\label{fig:MultistableChialvoone}
\end{figure}

In this section, we illustrate the \blue{phenomenon of coexisting attractors or multistability} in the case of \blu{the} Chialvo neuron under the action of electromagnetic flux. We observe \blu{the} coexistence of a chaotic attractor and periodic attractors with the change in electromagnetic flux $k$. Scanning the parameter space \blue{has} shown that \blue{different periodic solutions of high periods coexist} with the variation of the electromagnetic flux parameter $k$. Such coexistence is shown in Fig. \ref{fig:MultistableChialvoone}, where there is a coexistence of a chaotic attractor, a period-six solution, and \blue{a period-nine} solution. In \blue{Fig. \ref{fig:MultistableChialvoone} (a)}, a period-six solution is shown \blue{marked by yellow} triangles. In \blue{Fig. \ref{fig:MultistableChialvoone} (b)}, a period-nine solution is shown marked by red triangles. In \blue{Fig. \ref{fig:MultistableChialvoone} (c)}, a chaotic attractor is shown in black. In \blue{Fig. \ref{fig:MultistableChialvoone} (d)}, a slice of the basin of attraction of the attractors is shown at $\phi=0$ plane. To construct the basin of attraction plot, a fine grid of $x,y$ values was taken, and for each point sufficiently large number of iterations were taken and checked where it converges. \blue{If it converges to a period-six solution, it is marked in yellow, if to a period-nine solution, it is marked in red, if to a chaotic attractor, then is marked in black and finally if it diverges to infinity, it is marked in white.} Following such procedure produces Fig. \ref{fig:MultistableChialvoone} (d). \blue{It is seen that much of the basin of attraction is dominated by the prevalence of chaotic attractor. The size of the basin for a period-nine solution is less as compared to the size of the basin of the period-six solution. The system shows two symmetrical kind of basin of attraction in red, even though the map is asymmetric. Observe that at first, it might seem that the basin structure is regular, but there is a complicated structure associated \blu{with} it}. To illustrate the striking complexity of the basin of attraction structure, we next consider zoomed-in regions of the basin of attraction plot near the regions of period-six and period-nine solutions and have revealed \blue{that the basin structure is indeed highly intermingled.}

\begin{figure}[!htbp]
\begin{center}
\includegraphics[width=0.9\textwidth]{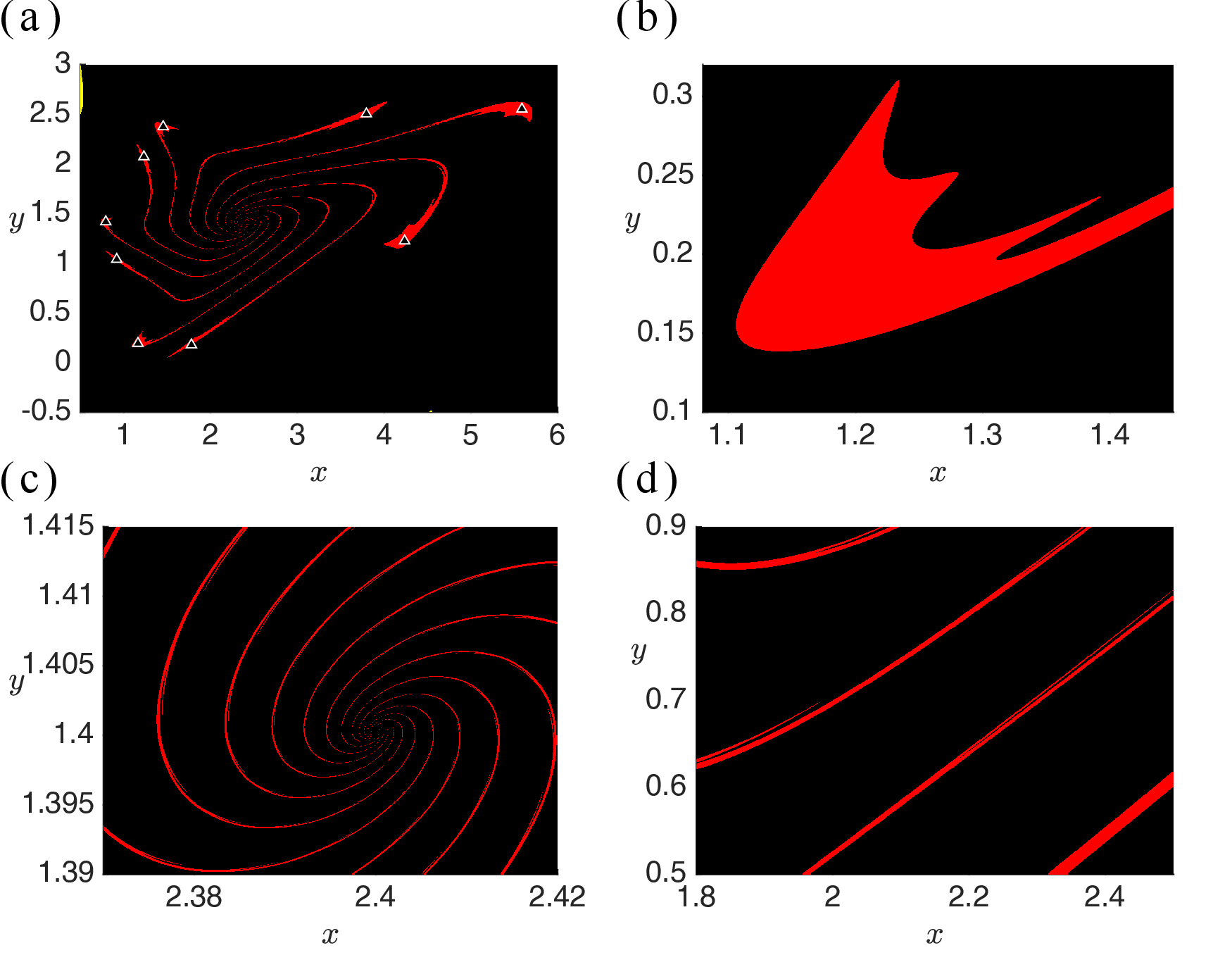}
\end{center}
\caption{Basin of attraction of the period-nine solution marked in red. The periodic points of period-nine are marked with black triangles. \blue{In (a), the period-nine solution is shown in triangles and they completely lie inside the basin. In (b), the shrimp structure associated with the basin structure is shown. In (c), the spiral structure of the basin is shown. It gets intermingled as one goes towards the center of the spiral structure. In (d), shrimp legs of the basin are shown. The parameters are set as $a=0.6,b=0.6,c=2,k_{0}=0.28,k_{1}=0.1,k_{2}=0.2,k=0.002,\alpha=0.1,\beta=0.2$.} }
\label{fig:BasinPeriodNine}
\end{figure}

A zoomed in version of the period-nine solution is shown in Fig. \ref{fig:BasinPeriodNine}. In \blue{Fig. \ref{fig:BasinPeriodNine} (a)}, \blue{the} periodic points are marked by black triangles. Observe that basin of attraction of the period-nine solution appears to be a spiral like structure. In \blue{Fig. \ref{fig:BasinPeriodNine} (b)}, a zoomed in version near one of the spiral arms is shown confirming the shrimp structure. In \blue{Fig. \ref{fig:BasinPeriodNine} (c)}, we zoom in at the center of the spiral region and we still get a spiral, confirming the spiral structure of the basins. In \blue{Fig. \ref{fig:BasinPeriodNine} (d)}, a zoomed in version is shown near the spiral arms showing that they are shrimp arms self repeating. 

\begin{figure}[t!]
\begin{center}
\includegraphics[width=0.7\textwidth]{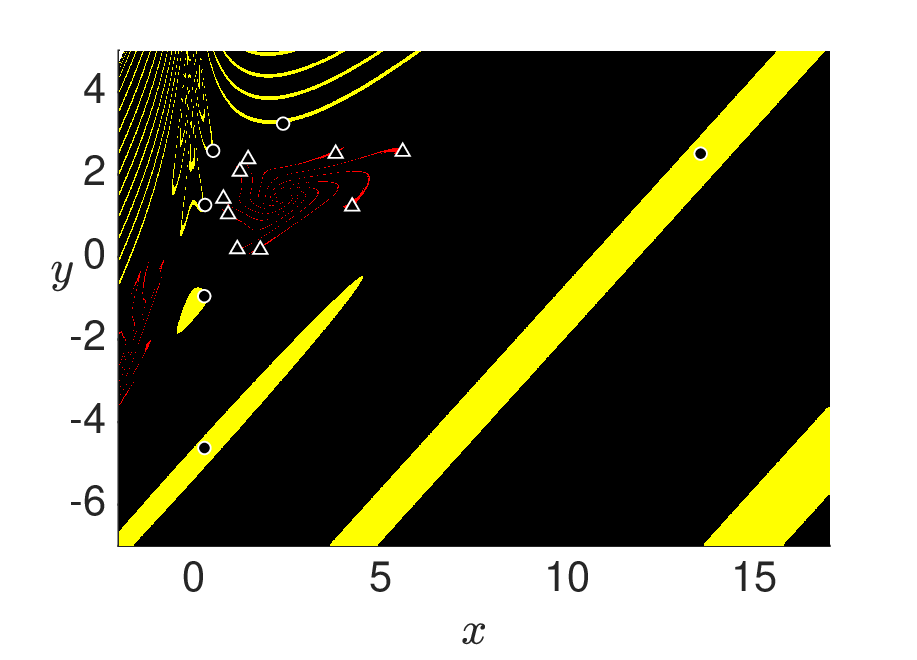}
\end{center}
\caption{Basin of attraction of the period-nine solution marked in red, \blue{period-six solution marked in yellow, and that of the chaotic attractor in black. Observe that the periodic points sit inside the basin of attraction confirming their asymptotic stability}. The periodic points of period-nine are marked with black triangles and the periodic points of the period-six are marked with black circles.  \blue{The parameters are set as $a=0.6,b=0.6,c=2,k_{0}=0.28,k_{1}=0.1,k_{2}=0.2,k=0.002,\alpha=0.1,\beta=0.2$.}}
\label{fig:BasinPeriodSix}
\end{figure}

The basin of attraction together with periodic points is shown in Fig. \ref{fig:BasinPeriodSix}. To avoid confusion, we have marked the period-nine solutions by triangles and period-six solutions by circles. Observe that the periodic points sit at the center of their respective  basins of  attraction confirming their asymptotic stability. \blue{Thus the period-nine and period-six solutions shown are both asymptotically stable}. We \blue{next show} the zoomed regions near the tip of the basin of attraction for $10<y<21$.

\begin{figure}[t!]
\begin{center}
\hspace*{-0.8cm}\includegraphics[width=1.1\textwidth]{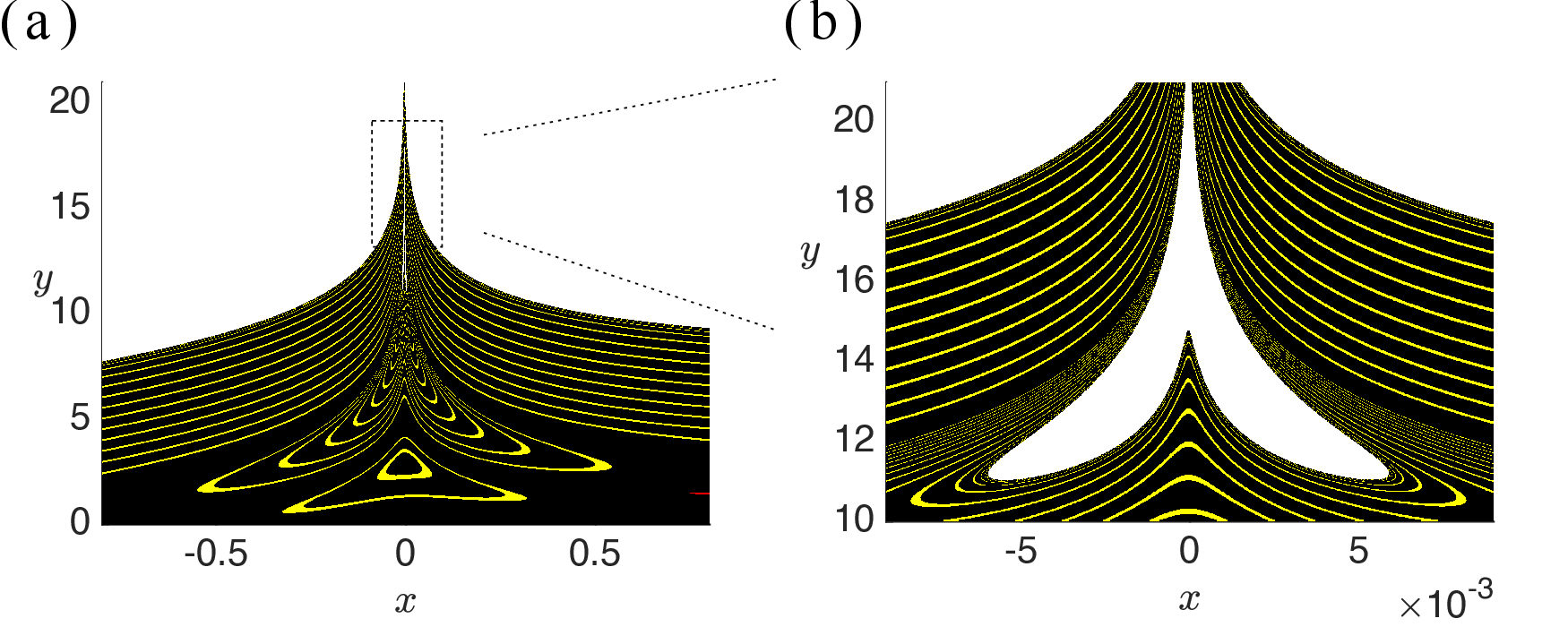}
\end{center}
\caption{A zoomed in version of a part of Fig. \ref{fig:MultistableChialvoone} above $y>0$ and $-0.8<x<0.8$. In (a), we can see the intermingled basins of chaotic attractor and period-six solutions, with a small region of divergence in the middle which is enlarged in (b). In (b), we observe a similar shape and also that the boundaries of the basins are complicated.  \blue{The parameters are set as $a=0.6,b=0.6,c=2,k_{0}=0.28,k_{1}=0.1,k_{2}=0.2,k=0.002,\alpha=0.1,\beta=0.2$.}}
\label{fig:CollageMultistable}
\end{figure}

In Fig. \ref{fig:CollageMultistable} (a), \blue{we} observe a closed disconnected sea of \blue{period-six} basins and they start to get distorted with increase in the value of $y$ near the tip. \blue{To} explore the region near the tip, around $10<y<21$, observe the fractal like structure between the basins of the chaotic attractor and period-six solution with a void of divergence. 

\section{Bifurcation structures and antimonotonicity}
\label{sec:bifstruc}
\begin{figure}[!htbp]
\begin{center}
\includegraphics[width=0.8\textwidth]{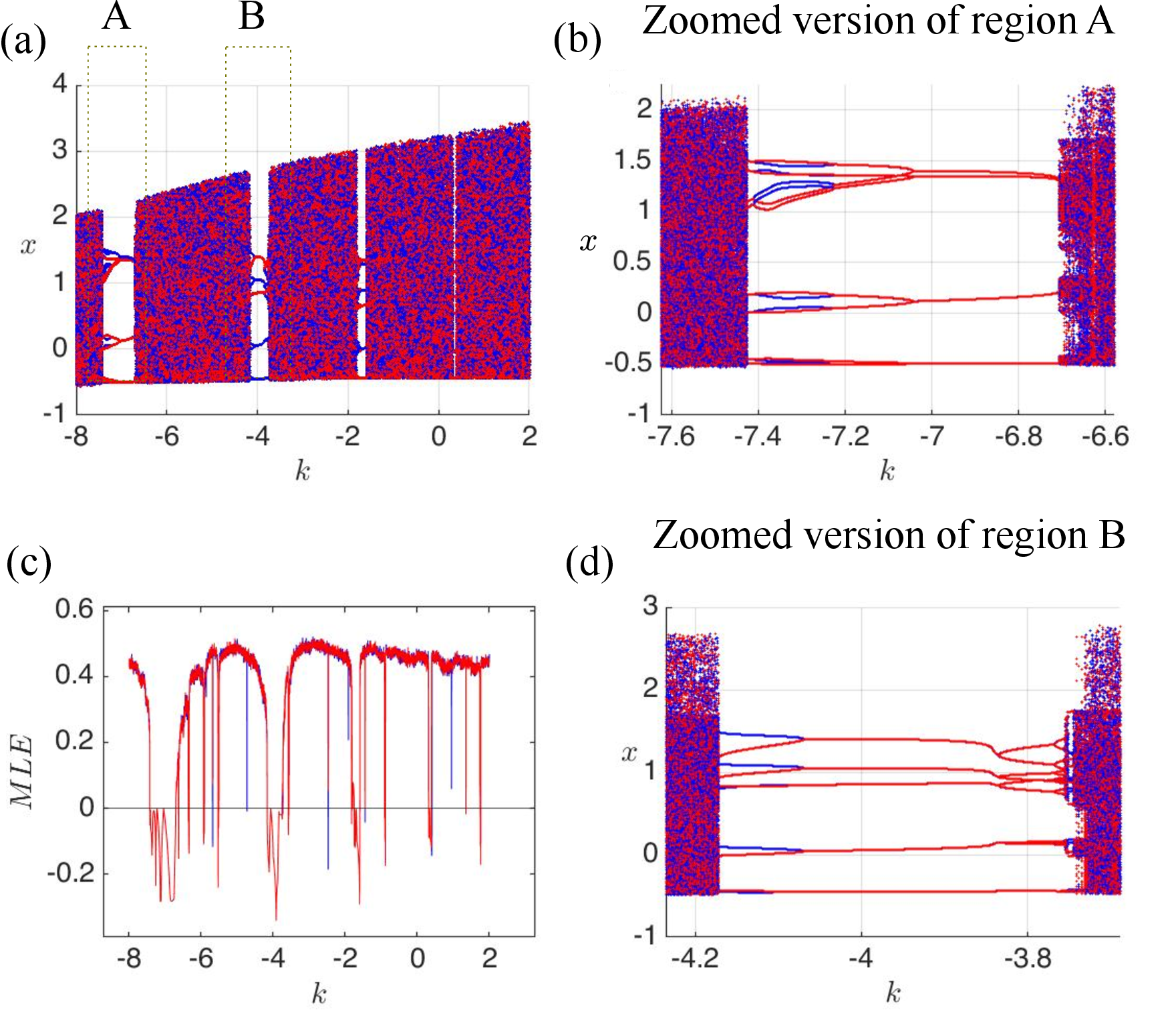}
\caption{Bifurcation diagram of $x$ with respect to $k$ in panel (a). A maximal Lyapunov exponent diagram is shown in panel \blue{(c)}. The blue colour denotes the forward bifurcation diagram and the red \blue{colour} denotes the backward bifurcation diagram. The parameter values are fixed as $a = 0.5, b=0.4, c = -0.89, k_{0} = -0.44, \alpha = 0.1, \beta = 0.1, k_{1} = 0.1, k_{2} = 0.2$. }
\label{fig:bifLyapK}
\end{center}
\end{figure}

In this section, we present a bifurcation analysis of the map \eqref{eq:ChialvoMag} with respect to electromagnetic flux parameter \blue{$k$}.  In Fig. \ref{fig:bifLyapK} (a), a one-parameter $x-k$ bifurcation diagram is constructed via both forward and backward continuation. In order to simulate the bifurcation diagram, for each value of $k$, \blue{some number of} final iterates \blue{(here last $100$ points were taken)} of the state variable $x$ are considered and plotted. A fine range of $k$ values \blu{was} taken in the interval of $k \in [-8,2]$. The forward continuation (marked in blue) is considered by increasing the value of $k$ from $-8$ to $2$ and the backward continuation (marked in red) is carried by decreasing the value of $k$ from $2$ to $-8$. Both the forward and backward continuation points were plotted in the same figure. This method provides with a benefit \blu{detecting} multistability. If both the forward and backward continuation points do not overlap completely, it suggests multistability in the system. \blue{Continuating} in different directions of parameter space reveals different \blu{long-term} behaviors or multistability. 

To explore the chaotic, periodic regions, we compute the Lyapunov exponents via $QR$-factorisation method \cite{GoFePa01, Muni22a}. In Fig. \ref{fig:bifLyapK}, we present the maximal Lyapunov exponent with respect to parameter $k$ for both the forward and backward continuation together.

In the region of $-5<k<-3$, we observe period-doubling and inverse period-doubling scenario. Through this bifurcation diagram, we can confirm the presence of periodic orbits of large periods. The period increases as $k$ increases. The bifurcation diagram in panel  \blue{Fig. \ref{fig:bifLyapK}}(a) is in correspondence with the Lyapunov spectrum diagram in panel  \blue{Fig. \ref{fig:bifLyapK}}(b).

In Fig. \ref{fig:bifLyapK} (b), we present a \blu{zoomed-in} version of region A marked in \blue{Fig. \ref{fig:bifLyapK}}(a), to showcase the periodic bubble in the bifurcation diagram. The formation of a bubble in the bifurcation diagram implies concurrent creation and destruction of periodic orbits. This phenomenon is known as antimonotonicity, motivated by the non-monotonic behavior of the bifurcation diagram \cite{Da92}. It is considered a fundamental behavior in nonlinear dynamical systems and is present in a large class of nonlinear systems. Experimental observation of antimonotonicity is observed in Chua's circuit \cite{Ko93}. The inevitable period reversal bubbles and forward period bubbles were reported.
In \blue{Fig. \ref{fig:bifLyapK}}(c), we observe periodic reversals in the bifurcation diagram by considering a zoomed in version of region $B$ in the bifurcation diagram. 

\begin{figure}[!htbp]
\begin{center}
\includegraphics[width=0.8\textwidth]{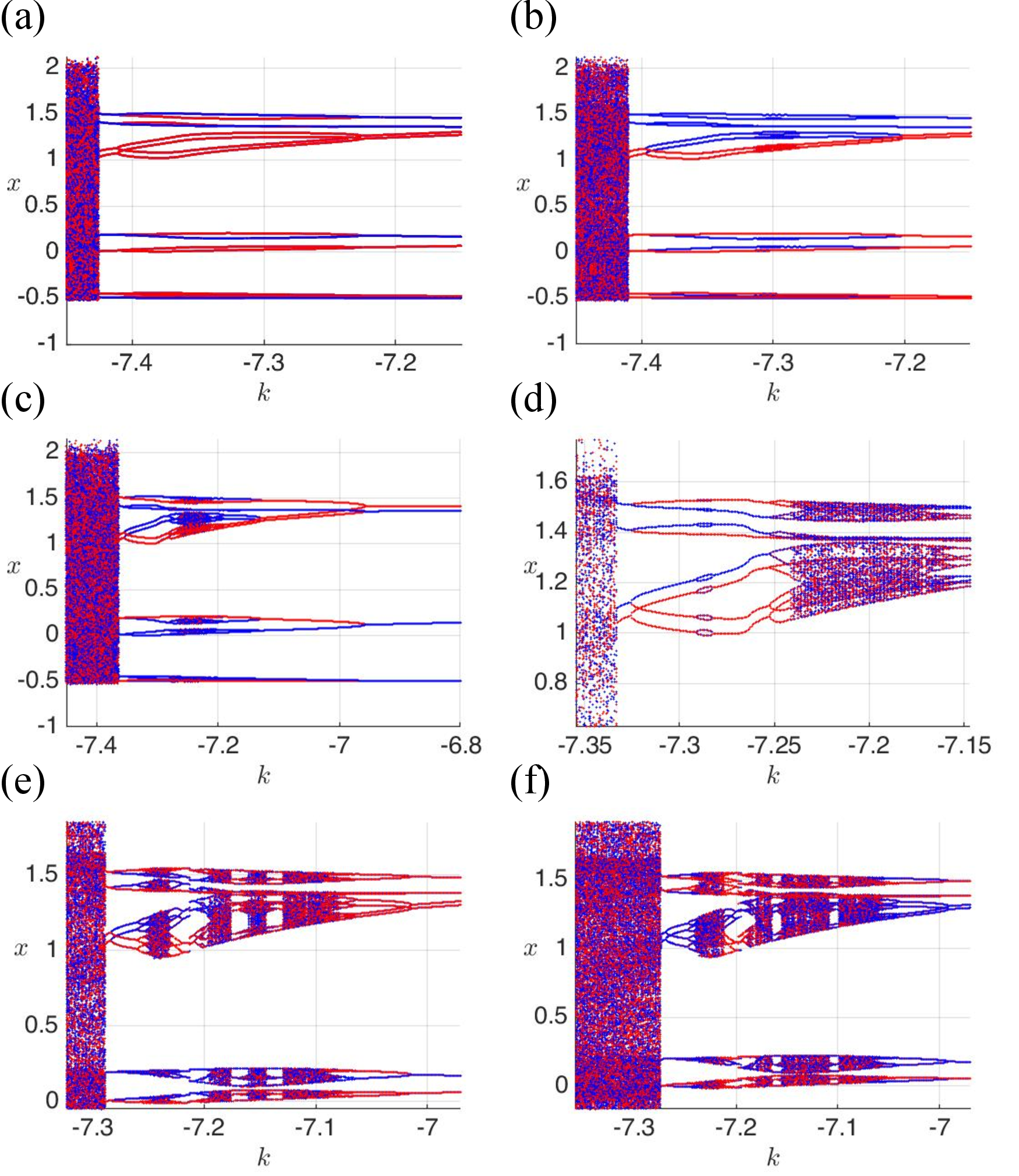}
\caption{Bifurcation diagram of $x$ with respect to $k$ in panel (a).  The blue colour denotes the forward bifurcation diagram and the red colur denotes the backward bifurcation diagram. The parameter values are fixed as $a = 0.5, b=0.4, c = -0.89, \alpha = 0.1, \beta = 0.1, k_{1} = 0.1, k_{2} = 0.2$. }
\label{fig:AntiMonBubble}
\end{center}
\end{figure}

For $-3.8 < k < -3.7$, we observe a period-doubling route to chaos.  We next explore the phenomenon of antimonotonicity route to chaos by the variation of parameter $k$. Such transition to chaos is shown in Fig. \ref{fig:AntiMonBubble}.  Zooming in the region in Fig. \ref{fig:bifLyapK} in the region of $-8<k<-6$, we observe antimonotonicity illustrated in Fig. \ref{fig:AntiMonBubble}. In \blue{Fig. \ref{fig:AntiMonBubble}}(a), we encounter the formation of periodic bubbles for $k_{0} = 0.44$. When $k_{0} = 0.441$, extra pair of bubbles are formed as shown in Fig. \ref{fig:AntiMonBubble}. When we increase $k_{0}$ slightly to $0.444$, we encounter the formation of chaotic bubbles in \blue{Fig. \ref{fig:AntiMonBubble}}(c). Further increase in $k_{0}$ to $0.446$, we encounter the presence of new periodic bubbles and chaotic bubbles in \blue{Fig. \ref{fig:AntiMonBubble}}(d). When $k_{0} = 0.449$ in \blue{Fig. \ref{fig:AntiMonBubble}}(e), a large portion of chaotic behavior appears with periodic reversals and many number of periodic bubbles are formed. Chaotic behavior with bubbles persists for $k_{0} = 0.45$. This illustrates the fact that reverse bubbles, antimonotonicity phenomenon occur infinitely many times near common parameter values \cite{Da92}. Rigorous mathematical analysis of antimonotonicity in nonlinear systems is yet to be uncovered.   

\begin{figure}[!htbp]
\begin{center}
\includegraphics[width=0.7\textwidth]{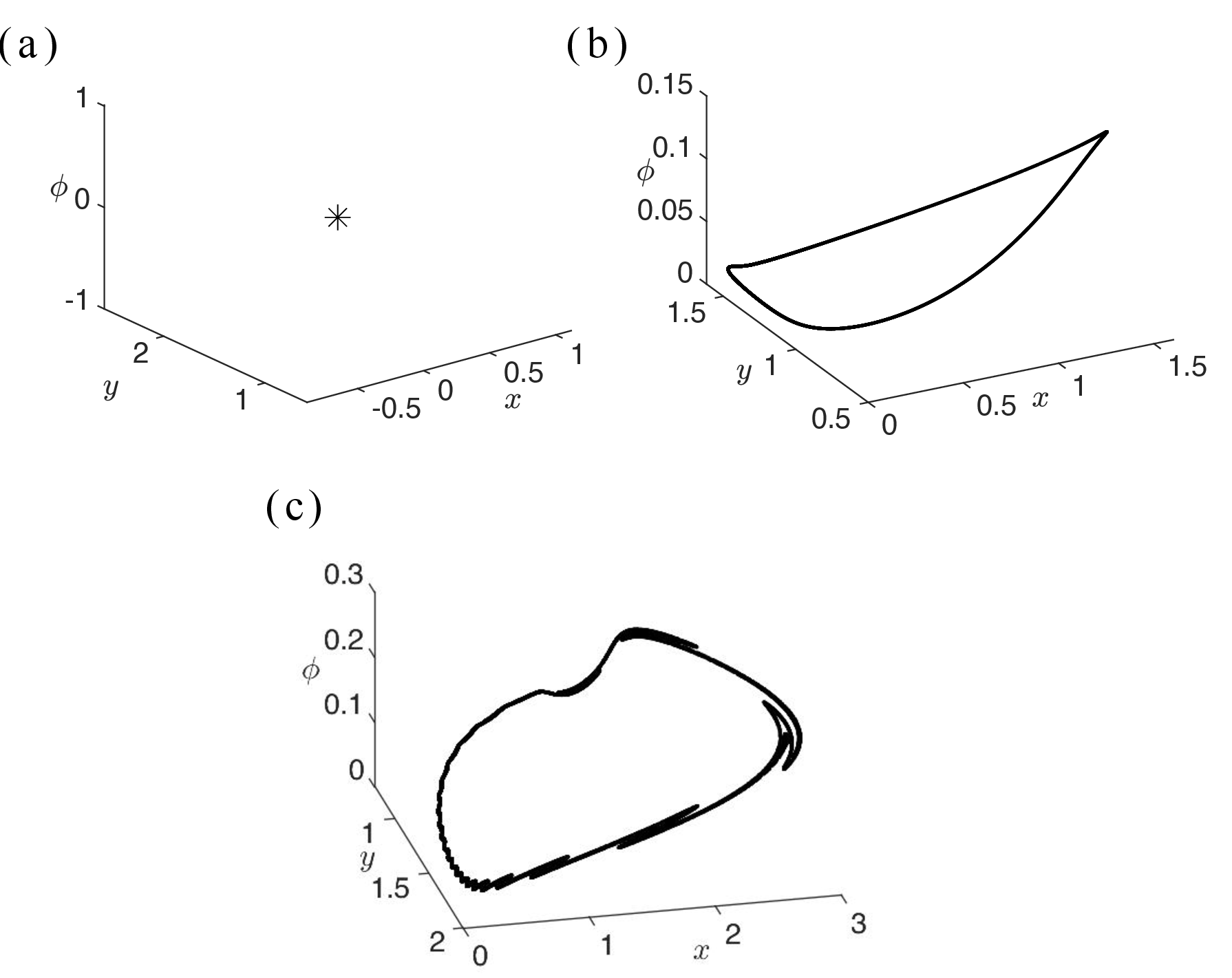}
\end{center}
\caption{In (a) a stable fixed point is shown in the $x-y-\phi$ \blue{phase} space for $a = 0.838$. After a supercritical Neimark-Sacker bifurcation, an attracting closed invariant curve is born as shown in (b) at $a=0.841$. A chaotic attractor is then formed when $a$ is increased to $0.88$. The parameters are set as $b=0.18, c = 0.28, k_{0}=0.06, k=-0.2, \alpha = 0.1, \beta = 0.2, k_{1} = 0.1, k_{2} =0.2$. }
\label{fig:ClosedInvariantCurves}
\end{figure}

There are some bifurcations \blu{that} are exhibited by higher dimensional maps and cannot be exhibited by the two-dimensional maps, for example, torus doubling bifurcations \cite{Kan84}. Here we showcase another route to chaos: attracting invariant closed curve to chaos \cite{Chr03}. We show that Chialvo neuron \eqref{eq:ChialvoMag} undergoes the route of invariant closed curve to chaos by varying the parameter $a$. In Fig. \ref{fig:ClosedInvariantCurves} (a), the Chialvo map \eqref{eq:ChialvoMag} has a stable fixed point when $a=0.838$. After a supercritical Neimark-Sacker bifurcation, an attracting closed invariant curve is born while the fixed point goes unstable when $a=0.841$. When $a$ is further increased to $0.88$, the attracting closed invariant curve is destroyed, and a chaotic attractor is born. 

\subsection{Evolution of fingered attractors}
In this section, we will explore various attractors observed in Chialvo map under the variation of electromagnetic flux parameter $k$. Specifically we showcase fingered attractors with different number of fingers as $k$ varies. Similar fingered attractors were found in mechanical impact system \cite{Seth20}, buck converters \cite{Foss96}. To our knowledge, such fingered attractors are reported first in this article for discrete neuron maps. 

\begin{figure}[!htbp]
\begin{center}
\includegraphics[width=0.9\textwidth]{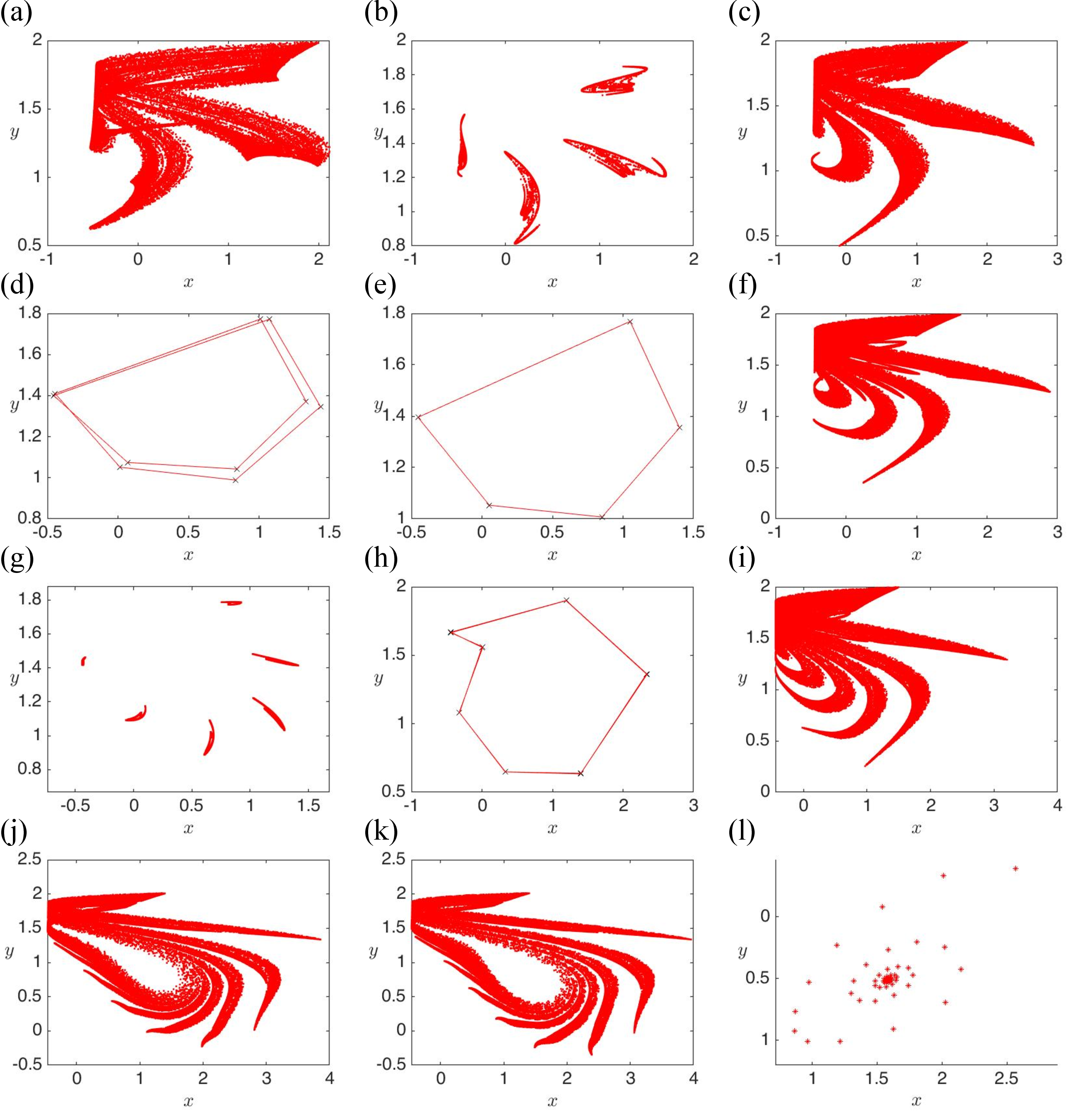}
\caption{Evolution of fingered attractors as electromagnetic flux parameter $k$ varies. The parameters are set as $a=0.5, b=0.4, c=0.89, k_{0}=-0.44, k_{1}=0.1, k_{2}=0.2, \alpha = 0.1, \beta = 0.1$.}
\label{fig:FingeredAttr}
\end{center}
\end{figure}

In Fig. \ref{fig:FingeredAttr} (a) for $k=-7.5$, a three fingered chaotic attractor is shown with three fingers. As $k$ increases to $-6.7$, four disjoint attractors are observed in (b). When $k=-4.4$, a four-fingered attractor is observed (c). Notice that the number of disjoint chaotic attractors in (b) is the same as the number of fingers in (c). 
The four fingered attractor is destroyed as $k$ is increased to $-4.1$ and a stable  double-round period-$10$ orbit is born as shown in Fig. \ref{fig:FingeredAttr} (d).  At $k=4$, a stable period-$5$ orbit was detected through period-halving bifurcation in (e). 
With further increase in the $k$ value to $-2.9$, we find the number of fingers in the chaotic attractor increases to five, see Fig. \ref{fig:FingeredAttr} (f). At $k=-1.8$, we observe six disjoint attractors as shown in Fig. \ref{fig:FingeredAttr} (g). A six-fingered chaotic attractor at $k=-1.2$ was detected in Fig. \ref{fig:FingeredAttr} (i). For $k=-1.7$, a stable periodic orbit of period $12$ is found. At $k=-1.6$, the period gets halved to six, a stable period-six orbit was found, Fig. \ref{fig:FingeredAttr} (h).  For $k=-0.9$, a stable period-$14$ orbit is detected. For $k=-0.3$, a seven fingered strange attractor is detected, see Fig. \ref{fig:FingeredAttr} (i). 

As we see, the number of fingers in the chaotic attractor increases with \blu{an} increase in the value of electromagnetic flux $k$. A period-$14$ is detected at $k=0.34$. It then gets halved at $k=0.45$, and hence a stable period-seven attractor is found. At positive values of $k$, for $k=5$, we observe an eight-fingered attractor, see Fig. \ref{fig:FingeredAttr} (j). At $k=6$, a ten fingered attractor is detected, see Fig. \ref{fig:FingeredAttr} (k). With a further increase in the value of $k > 6$, we observe that the system settles down to a stable fixed point and the membrane potential goes to a stable fixed point, see Fig. \ref{fig:FingeredAttr} (l).

A conjecture we formulate on the number of fingers in fingered attractors  based on the observation of the evolution of attractors in Fig. \ref{fig:FingeredAttr} is that the number of fingers in the fingered chaotic attractor is the same as the number of disjoint chaotic attractors which gets destroyed just before the formation of the chaotic fingered attractor. Moreover, much information about the shape of the fingered chaotic attractor can be obtained by the closure of the unstable manifold of a saddle fixed point of the Chialvo map \eqref{eq:ChialvoMag}.  

\section{Numerical bifurcation analysis of Chialvo map}
\label{sec:numbif}
In this section we investigate the influence of parameters on 
the dynamics of the map \eqref{eq:ChialvoMag} via one- and two-parameter bifurcation analysis. We consider $k$ and $c$ as the main bifurcation parameters. The bifurcation diagrams were computed numerically using {\sc MatContM} \cite{KuMe19}. Codimension-1 and codimension-2 bifurcation types presented are summarised in Table~\ref{Tab:bifpoints}.
\begin{table}[htbp]
\centering
\caption{Abbreviations of codimension-1 and codimension-2 bifurcations}
\label{Tab:bifpoints} 
\begin{tabular}{ |p{5.5cm}|p{1cm}||p{5cm}|p{1cm}|  }
 \hline
 \multicolumn{4}{|c|}{Codimension-1} \\
 \hline
 Saddle-node (fold) bifurcation  & LP    & Neimerk-Sacker bifurcation&  NS\\
 \hline
Period-doubling (flip) bifurcation&   PD  &   &\\
  \hline
 \multicolumn{4}{|c|}{Codimension-2} \\
 \hline
 Cusp    &CP & Chenciner &  CH\\
 Generalized flip&   GPD  & Fold-Flip &LPPD\\
  Flip-Neimark-Sacker& PDNS  &  Fold-Neimark-Sacker  &LPNS\\
 1:1 resonance& R1 & 1:2 resonance&R2\\
  1:3 resonance& R3  & 1:4 resonance&R4\\
 \hline
\end{tabular}
\end{table}

\begin{figure}[!htbp]
\begin{center}
\includegraphics[width=0.9\textwidth]{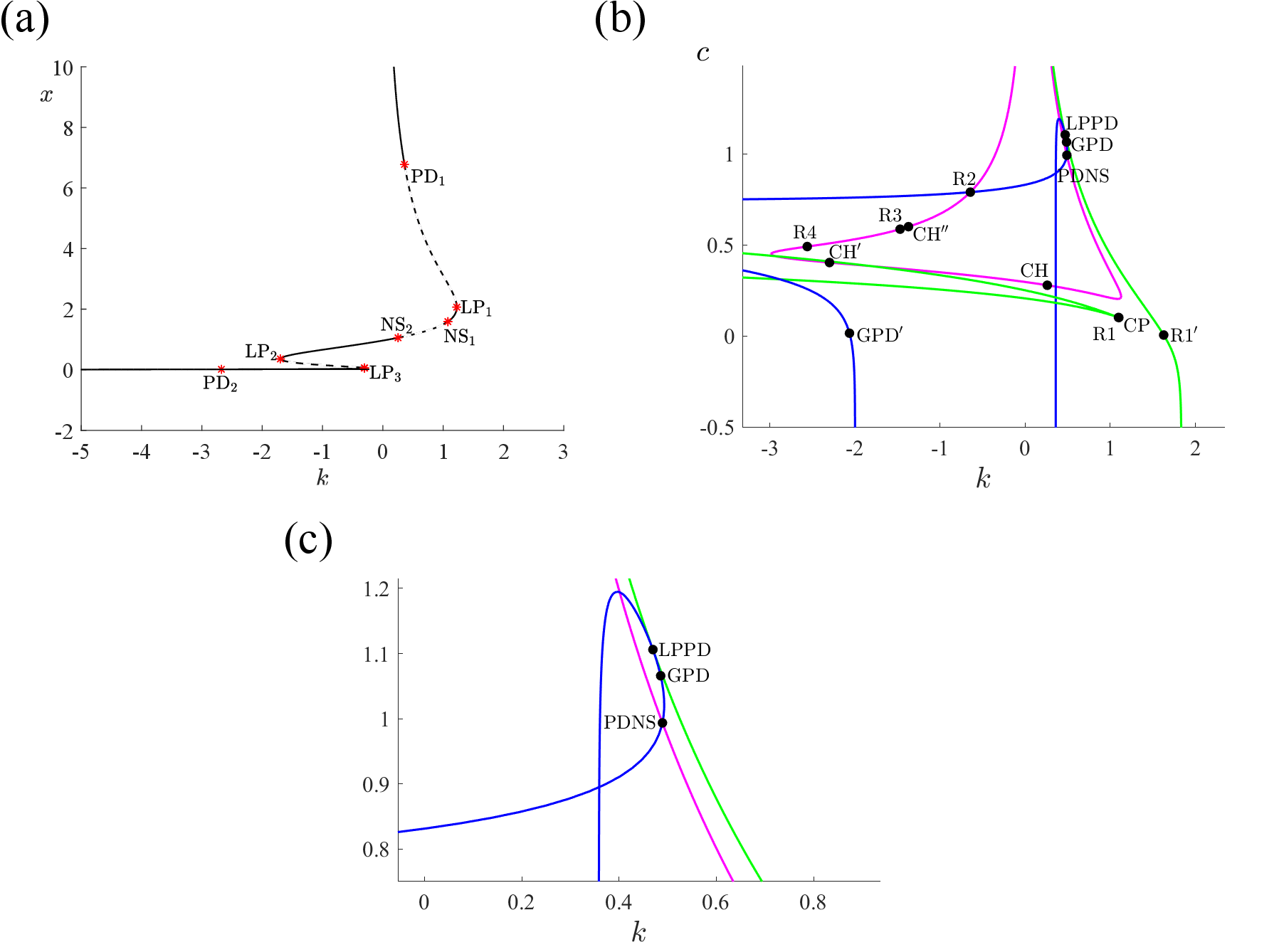}
\caption{(a) Codimension-1 bifurcation diagram of the map \eqref{eq:ChialvoMag} with $k$ as the  bifurcation parameter. Solid [dashed] black curves correspond
to stable [unstable] fixed points. (b) \blue{Codimension-2} bifurcation diagram in $(k,c)$-parameter plane. The green,blue, and magenta curves are the loci of the LP, PD, and NS bifurcations in (a), (b), and (c) respectively. The labels for the codimension-1 and -2 bifurcations are explained in Table~\ref{Tab:bifpoints}.}
\label{fig:NumBif1}
\end{center}
\end{figure}

Fig.~\ref{fig:NumBif1} (a) is a codimension-1 bifurcation diagram of $x$ as $k$ is varied with other parameters fixed. For sufficiently low values of $k$, the system has a single fixed point. As we increase the value of $k$, a subcritical period-doubling bifurcation ${\rm PD}_{2}$ with normal form coefficient $1.2164e^{+02}$ appears along the solution branch. Further increasing $k$, the fixed point loses stability in a saddle-node bifurcation ${\rm LP}_{3}$. The unstable branch emanates from the ${\rm LP}_{3}$ folds back at another saddle-node bifurcation $\rm {LP}_{2}$ to create a stable branch of fixed point, thus between the two saddle-node bifurcations $\rm {LP}_{2}$ and $\rm {LP}_{3}$ there are three fixed points: two stable and one unstable fixed points. Upon increasing the value of $k$, the stable fixed point in the upper branch loses stability in a supercritical Neirmark-Sacker bifurcation ${\rm NS}_{2}$ with normal form coefficient $-1.0239e^{-03}$. However, with further increasing $k$, a subcritical Neimark-Sacker bifurcation with normal form coefficient $6.5426e^{-02}$ appears, and further continuation results in saddle-node $\rm {LP}_{1}$ and period-doubling $\rm {PD}_{1}$ bifurcations along the solution curve.

Next we compute the two-parameter bifurcation analysis of the map \eqref{eq:ChialvoMag} by varying parameters $k$ and $c$. Fig.~\ref{fig:NumBif1} (b) is a bifurcation diagram in $(k,c)$-plane. The figure is composed of the curves of codimension-1 bifurcations in Fig.\ref{fig:NumBif1} (a). The blue, green and magenta curves are the loci of the period doubling bifurcation PD, saddle-node bifurcation LP and Neimark-Sacker bifurcation NS, respectively. 

For sufficiently large values of $c$, there exist Neimark-Sacker and saddle-node bifurcations, see Fig.~\ref{fig:NumBif1} (c) which shows a magnification of Fig.~\ref{fig:NumBif1} (b). As we decrease the value of $c$ a fold-flip bifurcation, denoted \blue{by} ${\rm LPPD}$, occurs on the saddle-node curve. At this codimension-2 point, the saddle-node LP curve collides tangentially with the period-doubling PD curve. Thus below the ${\rm LPPD}$, apart from the saddle-node and Neimark-Sacker bifurcations, there is period-doubling bifurcation. Observe also along the locus of the PD is a generalised period-doubling  bifurcation ${\rm GPD}$, at this codimension-2 point the locus of a  subcritical PD bifurcation meeets the locus of a supercritical PD bifurcation. Below the ${\rm GPD}$ point, appear 1:2 resonance ${\rm R2}$ and flip-\blue{Neimark}-Sacker ${\rm PDNS}$ bifurcations. The PD curve intersects the Neimark-Sacker curve at these codimension-2 points. As $c$ is decreased further, we observed 1:3 resonance ${\rm R3}$, 1:4 resonance ${\rm R4}$, and three Chenciner bifurcations, denoted by ${\rm CH}$, ${\rm CH}'$, and ${\rm CH}''$, appear on the Neimark-Sacker curve. The locus of a subcritical Neimark Sacker bifurcation meets the locus of a supercritical Neimark-Sacker bifurcation at the Chenciner bifurcation point. Also, apart from the existing loci of NS, PD, and LP bifurcations two additional LP curves appear. Finally, with $c$ decreasing further the two LP curves collide and annihilate in a cusp bifurcation point. Below \blue{these} points appear a generalised bifurcation  ${\rm GPD}'$ along the PD loci and two 1:1 resonance R1 and  ${\rm R1}'$ bifurcations.

\begin{figure}[!htbp]
\begin{center}
\includegraphics[width=0.9\textwidth]{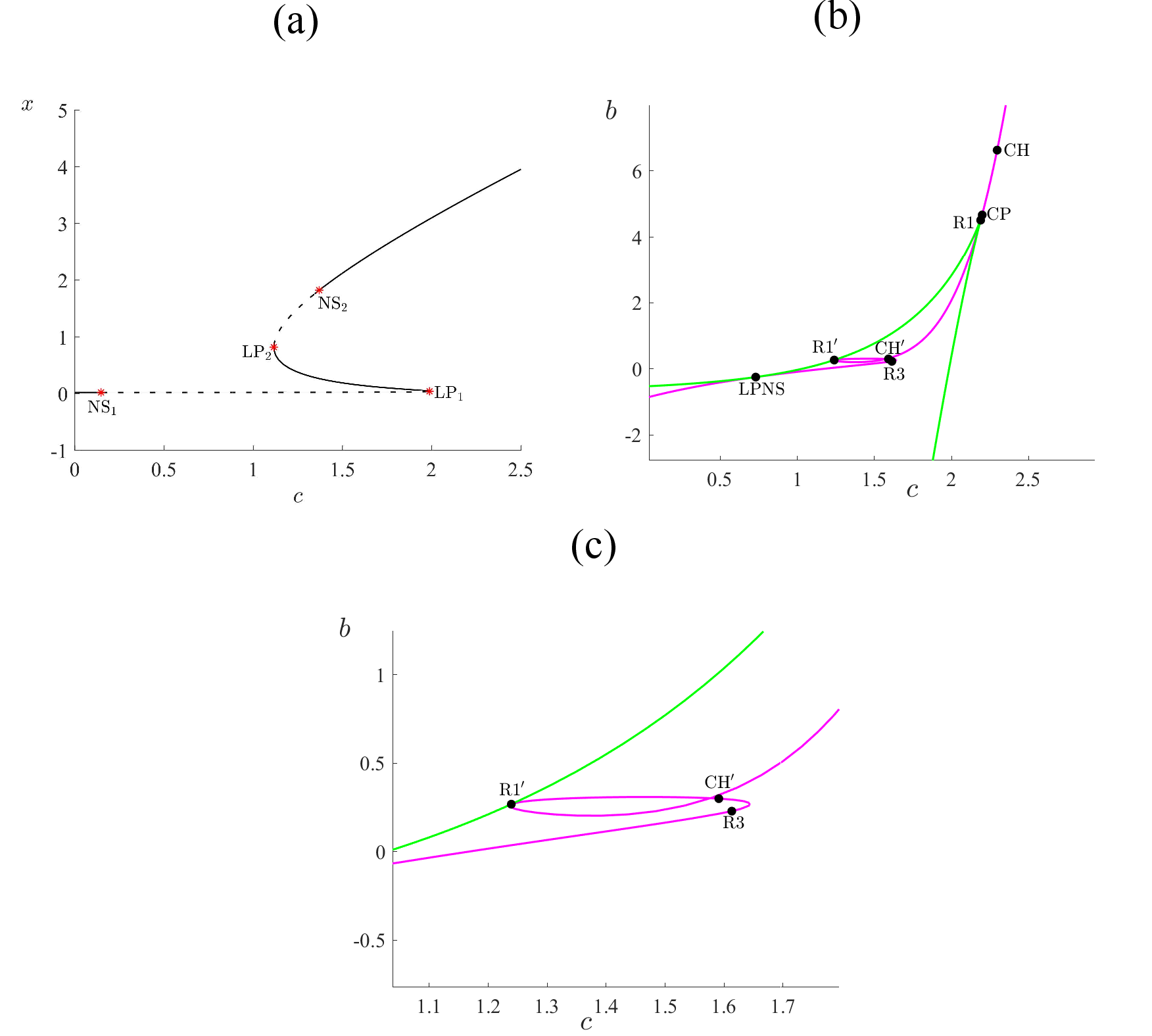}
\caption{(a) Codimension-1 bifurcation diagram of the map \eqref{eq:ChialvoMag} with $c$ as bifurcation parameter. Solid [dashed] black curves correspond
to stable [unstable] fixed points. (b) Codimension-2 bifurcation diagram in $(c,b)$-parameter plane. The green and magenta curves are the loci of the LP and NS bifurcations in (a), respectively. (c) The magnification of (b). The labels for the codimension-1 and -2 bifurcations are explained in Table~\ref{Tab:bifpoints}.}
\label{fig:NumBif2}
\end{center}
\end{figure}

Now we consider the bifurcation analysis of the map \eqref{eq:ChialvoMag} with $c$ as bifurcation parameter. A one-parameter bifurcation diagram of shown in Fig.~\ref{fig:NumBif2} (a). The map has a unique fixed point except between two
saddle-node bifurcations, ${\rm LP}_{1}$ and ${\rm LP}_{2}$ where three fixed points exist. To the left of ${\rm LP}_{1}$ the lower branch fixed point is unstable, and it \blu{gains} stability in a Neimark-Sacker bifurcation, denoted ${\rm NS}_{1}$. As we pass through the saddle-node bifurcation ${\rm LP}_{1}$ by increasing the value of $c$ the unstable lower branch folds to create a stable (middle) branch. The middle branch then folds back at the saddle-node bifurcation ${\rm LP}_{2}$ to produce an unstable fixed point which later \blue{regains} stability through a Neimark-Sacker bifurcation ${\rm NS}_{2}$. Thus, between ${\rm LP}_{1}$ and ${\rm NS}_{2}$, there is bistability, the coexistence of two stable fixed points.

A two-parameter bifurcation analysis of the map \eqref{eq:ChialvoMag} in $(c,b)$-parameter plane is shown in Fig.~\ref{fig:NumBif2} (b).  As we decrease the value of $b$ a Chenciner bifurcation $\rm CH$ appears on the NS curve. As the value of $b$ is decreased further, apart from the NS curve already observed, there are now two saddle-node curves, and the saddle-node curves collide and annihilate in a cusp bifurcation CP. Below the CP point is a 1:1 resonance bifurcation R1. With decreasing $b$, a Chenciner and a 1:3 resonance bifurcations, denoted by ${\rm CH}'$ and ${\rm R3}$ appear on the NS curve. Also, a 1:1 resonance bifurcation, denoted by ${\rm R1}'$, is observed and at this codimension-2 point the saddle-node curve coincides with the Neimerk-Sacker curve. Finally, as $b$ is decreased further appears a fold-Neimark-Saker bifurcation LPNS. At saddle-node and Neimark-Sacker curves \blu{intersect} at the LPNS point.

\begin{figure}[!htbp]
\begin{center}
\includegraphics[width=0.9\textwidth]{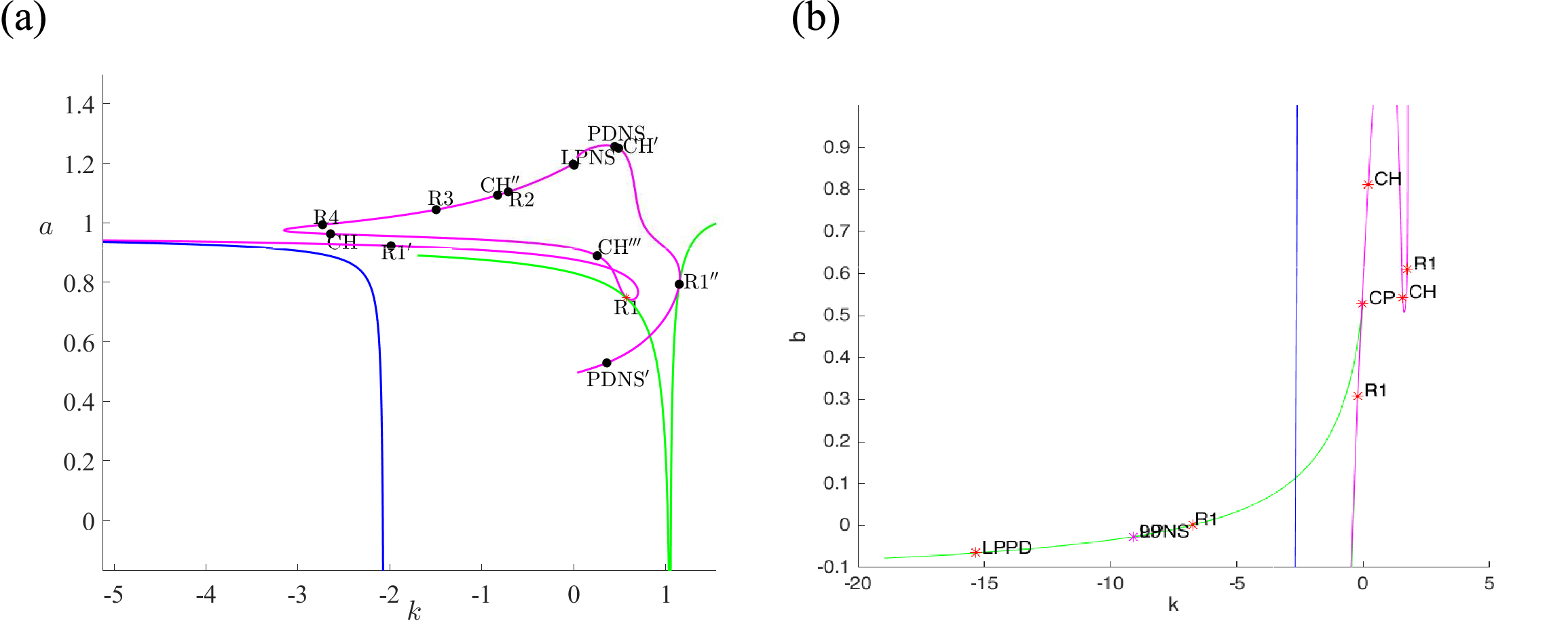}
\caption{\blue{Codimension-2 bifurcation diagram in (a) $(k,a)$-parameter plane with $b=0.05$; (b) $(k,b)$-parameter plane with $a=0.89$. The green, blue, and magenta curves are the loci of the LP, PD, and NS bifurcations, respectively. The parameters set as $c =0.28, k_{0}=0.03, \alpha = 0.5, \beta = 0.5, k_{1} = 0.2, k_{2} = 0.1$.}}
\label{fig:NumBif3}
\end{center}
\end{figure}

\blue{Lastly, we consider other parameter combinations. Figs.~\ref{fig:NumBif3} (a) and \ref{fig:NumBif3} (b) are codimension-2 bifurcation diagrams in $(k,a)$ and $(k,b)$ parameter spaces for the map \eqref{eq:ChialvoMag}. It is important to note that as the bifurcation parameters are varied simultaneously apart from the structure of the bifurcation curves similar codimension-2 bifurcations are observed as in Figs.~\ref{fig:NumBif1} (b) and \ref{fig:NumBif2} (b).
}

\section{Bursting and spiking features}
\label{sec:firing}
\begin{figure}[!htbp]
\hspace*{0.5cm}\includegraphics[width=0.7\textwidth]{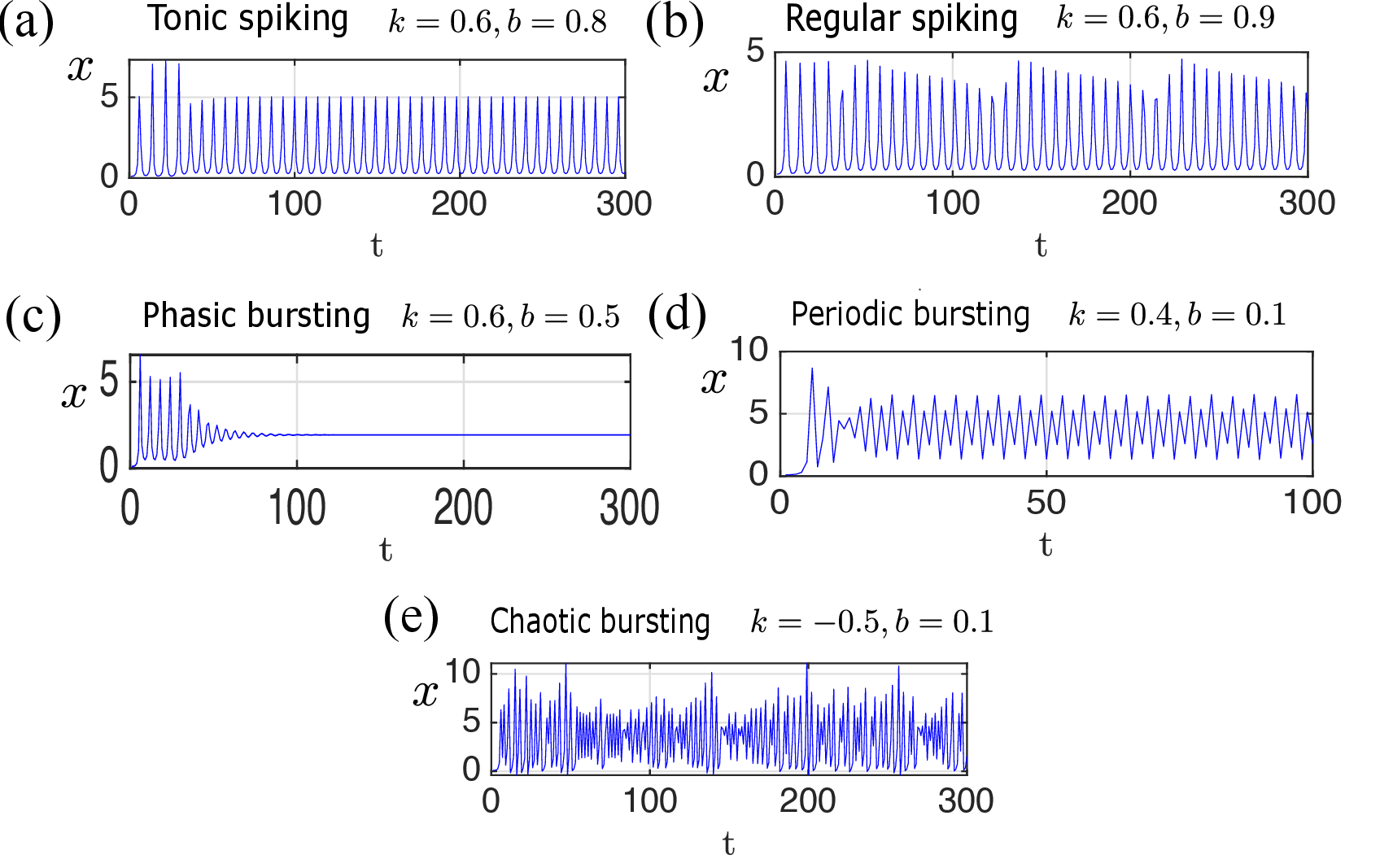}
\caption{Various spiking and bursting patterns exhibited by the Chialvo neuron model under the action of electromagnetic flux. In (a), tonic spiking is observed with intial high spike and then uniform spiking. In (b), regular spiking behavior is observed. In (c), Initial bursting followed by steady state depicting phasic bursting scenario is observed. In (d), periodic bursting is observed. In (e), chaotic spiking is shown.  The parameters set as $a=0.6, c = 1.4, k_{0}=0.1, \alpha = 0.1, \beta = 0.1, k_{1} = 0.1, k_{2} = 0.2$.}
\label{fig:BurstSpikeChialvo}
\end{figure}

 Bursting and spiking patterns \blu{is} an important behavior in the neurons.  It was shown in \cite{Muni22a} that \blu{discretized} Izhikevich neuron could exhibit many different spiking and bursting patterns by tuning the electromagnetic flux $k$. After performing bifurcation analysis, \blu{and} various routes to chaos, we explore the firing patterns exhibited by the Chialvo map \eqref{eq:ChialvoMag} under the action of electromagnetic flux $k$. In Fig. \ref{fig:BurstSpikeChialvo} (a), When $k=0.6, b=0.8$, tonic spiking is observed.

 When parameter $b$ is increased to $0.9$, a regular spiking pattern is observed in (b). When $b$ is instead decreased, neuron bursts initially for short time and then goes to steady state exhibiting phasic bursting behavior. In (d), a period-two bursting is observed when parameter $b$ is further decreased to $0.1$. In the \blu{negative-flux} region, specifically for $k=-0.5$, we observe a chaotic bursting pattern as shown in (e). 

\section{Chimera states in ring-star network of Chialvo map}
\label{sec:chimera}
\begin{figure}[htbp]
\hspace*{4cm}\includegraphics[width=0.5\textwidth]{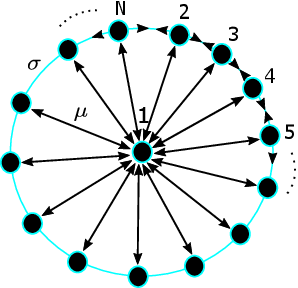}
\caption{Ring-star network of Chialvo neuron system. Each node $i=1 \ldots N$ represent a Chialvo neuron. The double-sided arrows signifies that the coupling is bidirectional in both ring and star connection. The coupling strengths of star and ring network are denoted by $\mu$ and $\sigma$ respectively.}
\label{fig:ChialvoRingStarNet}
\end{figure}

After the analyses are done in earlier sections on a single Chialvo neuron map, we will next focus our attention on a network of Chialvo neurons in this section. The network that we consider here has a ring-star configuration, first introduced in \cite{Muni20} and is illustrated in Fig.~(\ref{fig:ChialvoRingStarNet}). This configuration provides us with \blu{the} advantage of getting ring and star networks separately too. The mathematical model for the ring-star connected Chialvo neuron map under electromagnetic flux is defined as:
\begin{align}
    \label{eq:ChialvoDiscrete}
    x_m(n+1)&= x_m(n)^2e^{y_m(n)-x_m(n)} + k_0 + kx_m(n)M(\phi_m(n)) + \blue{\mu(x_m(n)-x_1(n))} + \\ \nonumber
            &+\frac{\sigma}{2R}\sum_{i=m-R}^{m+R}(x_i(n)-x_m(n)), \\ \nonumber
    y_m(n+1)&=ay_m(n)-bx_m(n)+c, \\ \nonumber
    \phi_m(n+1) &= k_1x_{m}(n)-k_2\phi_m(n),
\end{align}
whose central node is further defined as
\begin{align}
    \label{eq:ChialvoDiscreteCentral}
    x_1(n+1) &= x_1(n)^2e^{(y_1(n)-x_1(n))}+k_0+kx_1(n)M(\phi_1(n)) \\ \nonumber
            &+ \mu \sum_{i=1}^N(x_i(n)-x_1(n)), \\ \nonumber
    y_1(n+1) &= ay_1(n) - bx_1(n)+c, \\ \nonumber
    \phi_1(n+1) &= k_1x_1(n) - k_2\phi_1(n), 
\end{align}
having the following boundary conditions:
\begin{align}
    \label{eq:BoundaryConditions}
    \blue{x_{m+N}(n)= x_n(m)}, \\ \nonumber
    \blue{y_{m+N}(n) = y_n(m)}, \\ \nonumber
    \blue{\phi_{m+N}(n) = \phi_n(m).}
\end{align}

In (\ref{eq:ChialvoDiscrete}), $\sigma$ denotes the ring coupling strength between the neurons in the ring network and in (\ref{eq:ChialvoDiscreteCentral}), $\mu$ represents the coupling strength between the neurons and the central node of the star network. $R$ indicates the number of neighbors, and $M(\phi_m(n)) = \alpha+3\beta\phi_m(n)^2$ is the cubic memductance function. We have defined the size of the network to vary from $m=1$ to $m=N$ and simulate the results considering $N=100$. The parameters are set as $a = 0.89, b = 0.6, c = 0.28, k_0 = 0.04, k_1 = 0.1, k_2 = 0.2, \alpha = 0.1$, and $\beta = 0.2$.

\subsection{Ring network }
\begin{figure}[!htbp]
\begin{center}
\includegraphics[width=0.9\textwidth]{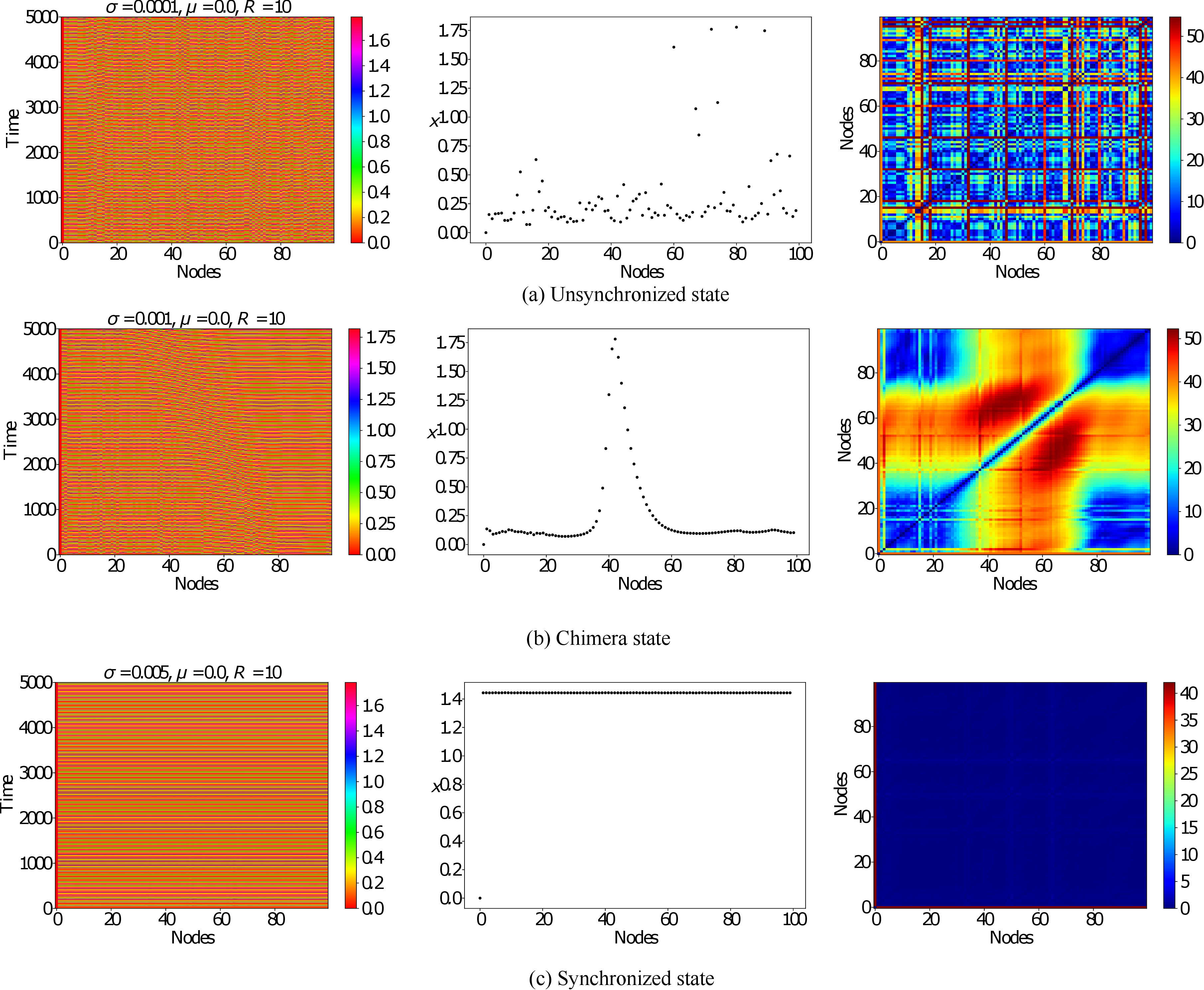}
\caption{Ring network of Chialvo neuron ($\mu = 0, \sigma > 0, k= 3.5$). In (a), unsynchronized behavior of the neurons is shown when $\sigma = 0.0001$. In (b), when $\sigma$ is increased to $0.001$, chimera state is observed. When $\sigma$ is further increased to $0.005$, \blue{synchronous} behavior is observed in (c).}
\label{fig:CombinedOneRing}
\end{center}
\end{figure}

Here, we \blu{consider} a simple ring network by considering $\mu=0$, i.e., eliminating the central node. We set $\sigma$ as the control variable and generate the figures representing the spatiotemporal dynamics of the network varying $\sigma$ from $0.0001$ to \blue{$0.005$}. We further set the coupling range to be $R = 10$ for our analysis, see Fig.~ \ref{fig:CombinedOneRing}.

We have used random initial conditions to make the system more complex. The leftmost plot gives the spatiotemporal dynamics, the middle plot the end state values of the nodes, and the rightmost plot the recurrence plot for the nodes.  When $\sigma = 0.0001$ in Fig. \ref{fig:CombinedOneRing}(a), unsynchronous pattern is formed evident from the end state values of the nodes and recurrence plot. When $\sigma $ is increased to $0.001$, we observe a \blu{chimera-like} state as there is a coexistence of synchronous and asynchronous nodes in the network, see Fig. \ref{fig:CombinedOneRing}(b). As we keep increasing $\sigma$, keeping $\mu = 0$, the system moves towards a synchronized state. Specifically when $\sigma = 0.005$, synchronized pattern is observed in Fig. \ref{fig:CombinedOneRing}(c). 


\subsection{Ring-star network }
\begin{figure}[!htbp]
\begin{center}
\includegraphics[width=0.9\textwidth]{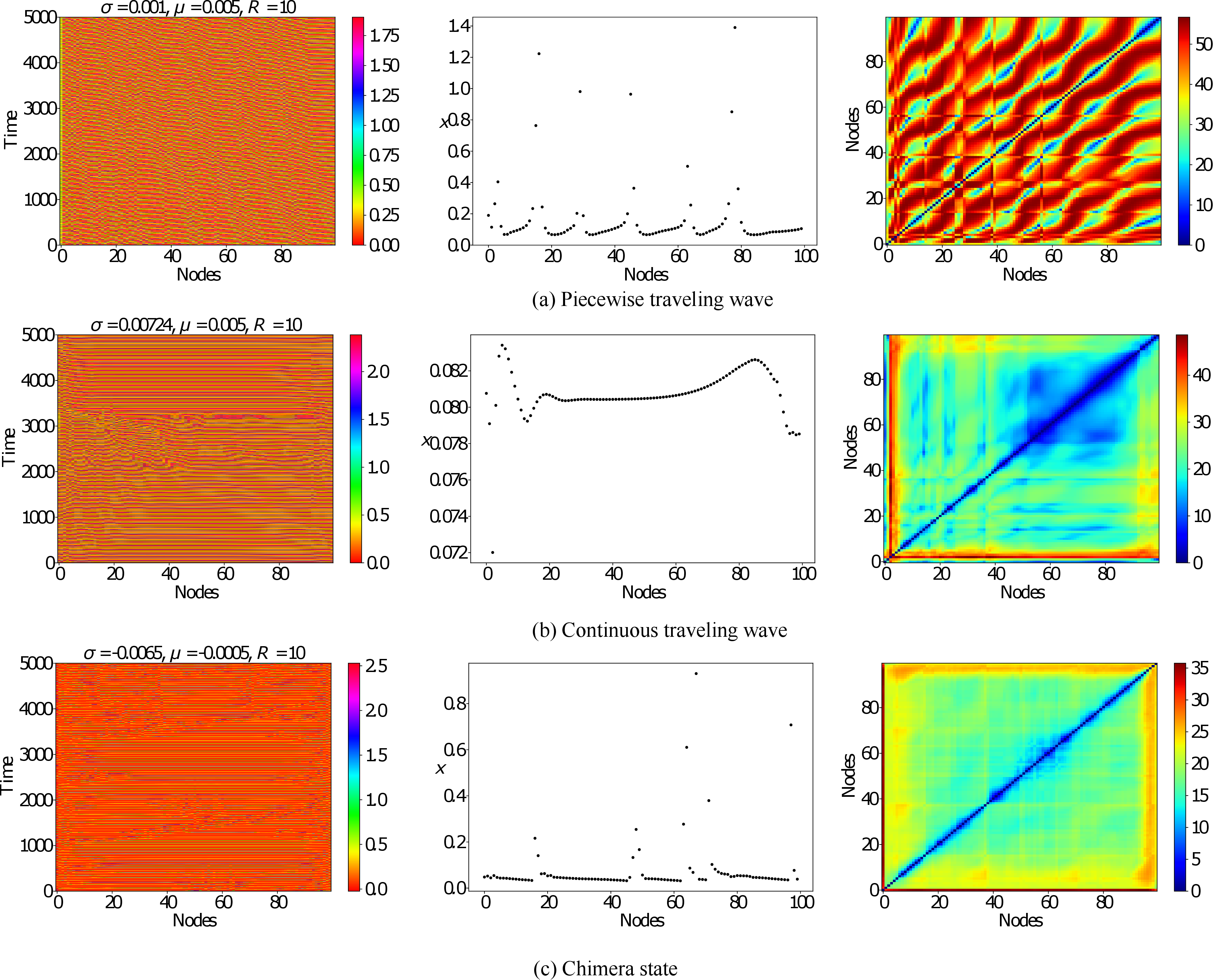}
\caption{Ring-star network of Chialvo neuron ($\mu > 0, \sigma \geq 0$) or ($\mu < 0, \sigma \leq 0$). In (a), there exists a piecewise wavy structure of the nodes for $\sigma = 0.001$. In (b), with increase in $\sigma$ to $0.00724$, we observe a continuous wavy structure of the nodes. In (a) and (b), we have $\mu>0, \sigma \geq0$ and $ k=3.5$. In (c), a chimera pattern of the nodes is shown where $\mu<0, \sigma = -0.0065$ and $k=-1$.}
\label{fig:CombinedOneRingStar}
\end{center}
\end{figure}

In this case, a ring-star network of \blu{the} Chialvo neuron model is considered. The coupling range is fixed at $R=10$. The  ring and star coupling strengths are varied to account for different spatiotemporal patterns. Various simulation runs have shown \blu{asynchronous} behavior  for low coupling strengths and synchronous behavior for larger coupling strengths.  Here we report the coexistence of synchronous and \blu{asynchronous} nodes (chimera state ) in Fig. \ref{fig:CombinedOneRingStar} (c). Next, we report two kinds of spatiotemporal structures which we observed and have classified them as a) piecewise wavy structures b) continuous wavy structures. In Fig. \ref{fig:CombinedOneRingStar}(a), a piecewise wavy \blue{structure} is revealed which is very near to be classified as a chimera state. In Fig. \ref{fig:CombinedOneRingStar} (b), a continuous wavy structure is evident from the end node plot.

\subsection{Star network }
\begin{figure}[!htbp]
\begin{center}
\includegraphics[width=0.9\textwidth]{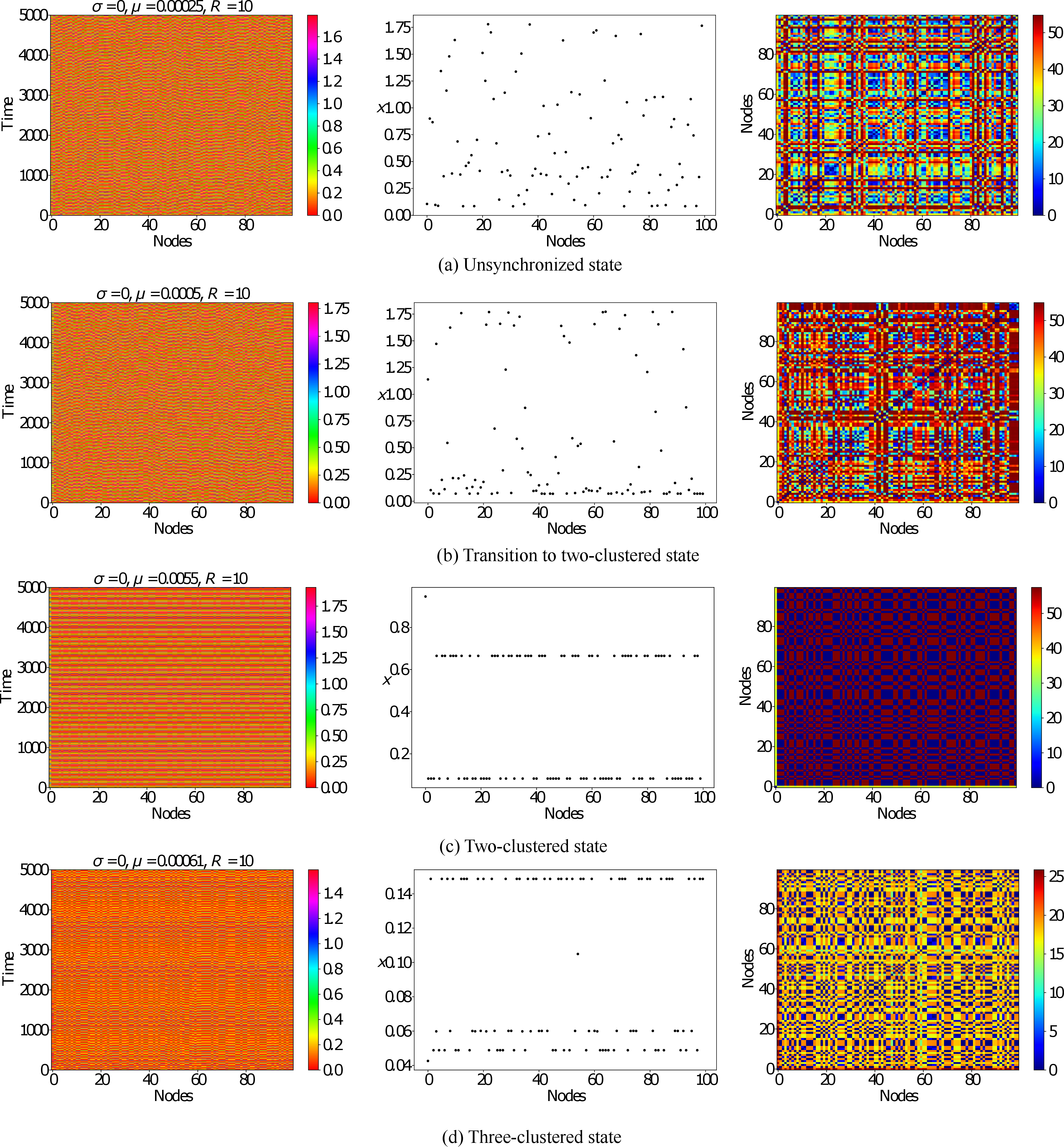}
\caption{\blue{Star network ($\mu > 0, \sigma = 0, k=3.5$ or $k = -1$). In (a), we observe asynchronous behavior of the nodes where $\mu =0.00025$ and $k=3.5$. In (b), we notice the transition to a two-clustered state with $\mu = 0.0005$.  In (c), we notice a perfectly two-clustered state with $\mu = 0.0055$. Finally in (d), we notice a perfectly three-clustered state where $\mu = 0.00061$ and $k=-1$}}
\label{fig:CombinedOneStar}
\end{center}
\end{figure}

In this section, we set the coupling strength between the nodes of the ring to $\sigma=0$, reducing our network to a star-network. Here $\mu$ is the control parameter.  Like the above two cases, we again fix \blu{the} coupling range $R=10$. For the case $\mu = 0.00025$, we record a highly unsynchronized behavior from the nodes (see Fig. \ref{fig:CombinedOneStar}(a)). For the case $\mu = 0.0005$, we recognize unsynchronized behavior transitioning towards a two-cluster state (see Fig. \ref{fig:CombinedOneStar}(b)) and finally for the case $\mu = 0.0055$, \blue{we observe the nodes having settled to a synchronized two-cluster state} in Fig. \ref{fig:CombinedOneStar} (c). The very small structures in the recurrence plot of Fig. \ref{fig:CombinedOneStar}(c) confirm the clusters of the two-clustered synchronized nodes.  We have noticed that with the change in the values of $\mu$ their is a merge of the nearest larger cluster forming the synchronized two-clustered state. 

\begin{figure}[!htbp]
\begin{center}
\includegraphics[width=0.5\textwidth]{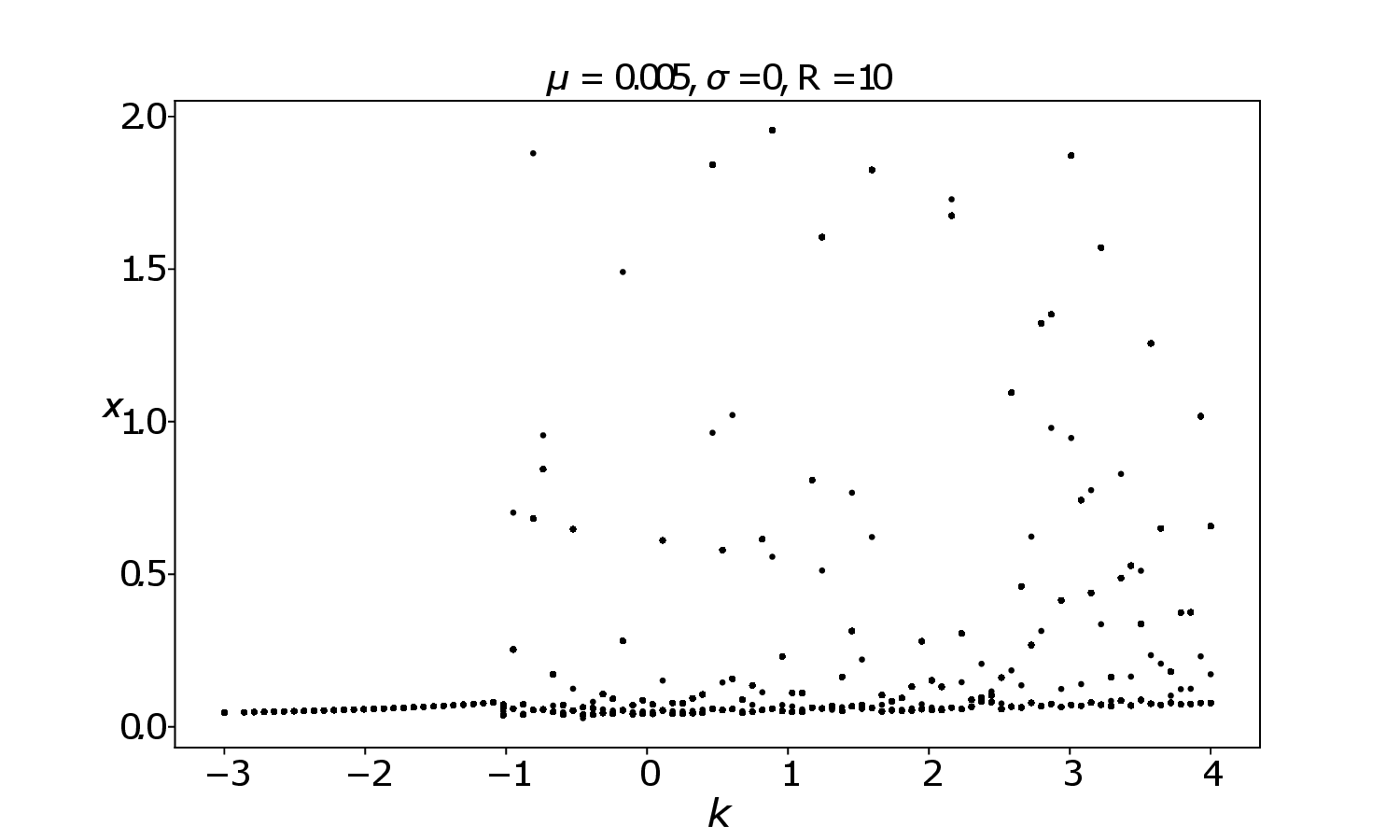}
\end{center}
\caption{$X-k$ plot for star network with $\mu = 0.005$ and $R=10$. This shows the intervals of $k$ for the prevalence of single-cluster, \blue{and} multi-cluster states. }
\label{fig:Xk_star}
\end{figure}

Synchronized cluster states are common in star networks \cite{MuPa18}. Can there be more than two cluster states in the star network of Chialvo neurons? \blu{To} investigate this further, we did a plot of $x$ state of all the nodes with the variation of electromagnetic flux $k$, see Fig. \ref{fig:Xk_star}. Observe that the system shows single cluster states for negative $k$ ($-3 < k < -1.3$). But as $k \geq -1$, there is a prevalence of more cluster states, as evident from the number of dots.  \blue{Such a case of three clustered state is shown in Fig. \ref{fig:CombinedOneStar} (d). Observe that the squares have reduced in size in the recurrence plot, and the number of squares \blu{has} instead increased in Fig. \ref{fig:CombinedOneStar} (d) when compared to Fig. \ref{fig:CombinedOneStar} (c).}

\subsection{Spatiotemporal patterns for negative coupling strengths}
After investigating ring-star, ring, \blu{and} star networks in the above sections, we present here a collection of piecewise and continuous wavy patterns we noticed while exploring spatiotemporal patterns in the regime of negative coupling strengths. Many simulations in the regime of both positive and negative strengths have shown that rich variety of spatiotemporal patterns are present in the case of negative coupling strengths. 

\blue{Negative (inhibitory) coupling strengths account for a significant proportion of neuronal connectivity in the human nervous system. The authors in \cite{TsKoKaPr18, TsHiHoPr16} have mentioned them in their course of simulations of the leaky Integrate-and-Fire (LIF) model. Also, negative coupling strengths are included via a rotational coupling matrix (See Eq. (2) in \cite{OmOmHoSc}) during simulations of FitzHugh Nagumo neuronal models.}

\begin{figure}[!htbp]
\begin{center}
\includegraphics[width=0.9\textwidth]{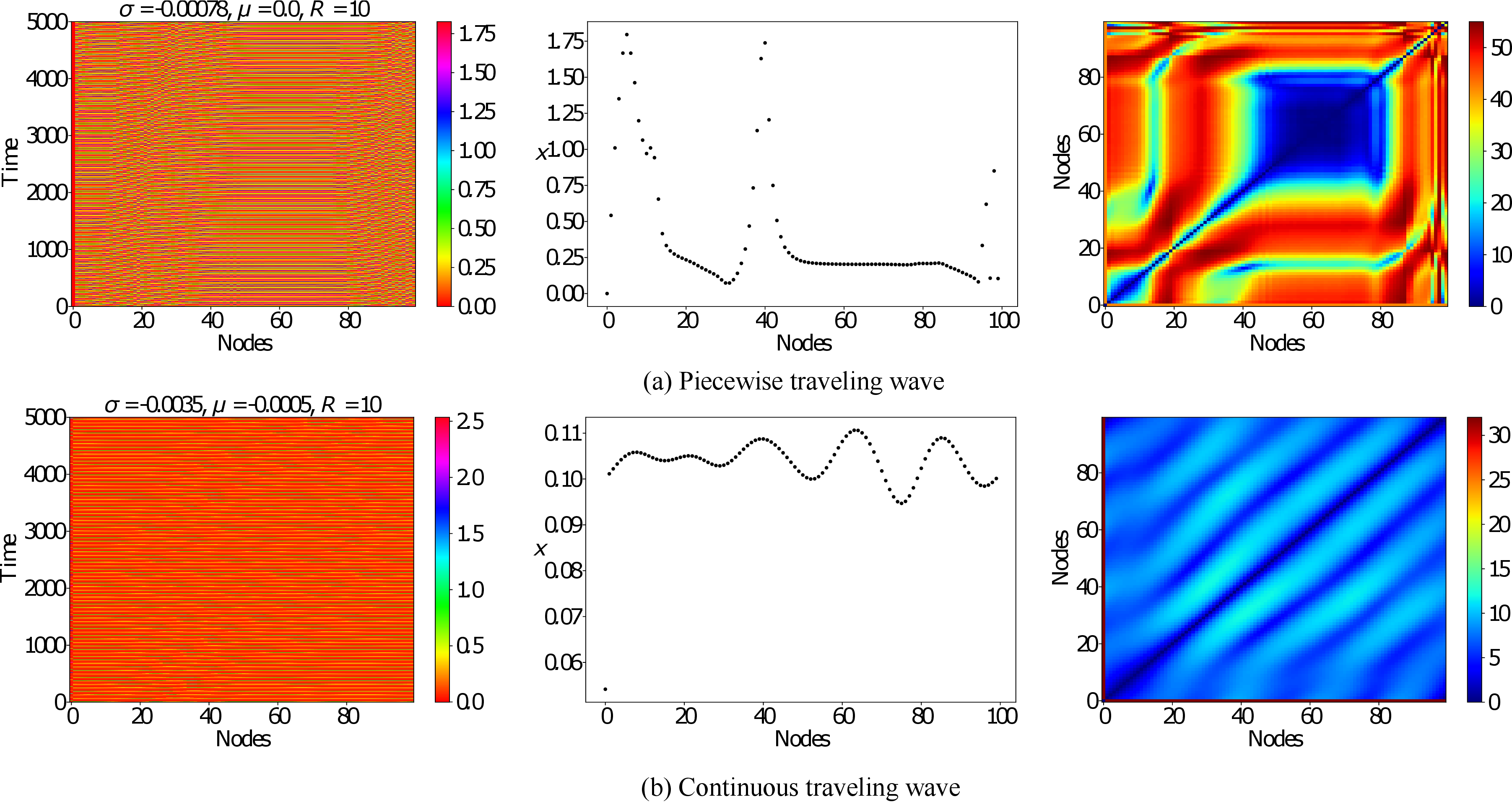}
\caption{{\color{black} Negative-strength patterns of Chialvo neuron ($\mu \leq 0, \sigma < 0$). In (a) we notice a piecewise-wavy pattern of the ring network nodes where $\mu = 0, \sigma < 0, k=3.5$. In (b) we observe a continuous-wavy structure of the ring-star network nodes with negative strengths where $\sigma<0, \mu<0, k=-1$.}}
\label{fig:NegativeStrengthsPattern}
\end{center}
\end{figure}

In Fig. \ref{fig:NegativeStrengthsPattern} (a), we account for a piecewise wavy structure when $\sigma = -0.00078$ and $k=3.5$. Also observe that in this state the recurrence plot is quite complex as well and topologically different from the cases of chimera state and unsynchronized state. At first glance, it is tempting to classify it as an unsynchronized behavior but since this appears to have an order to it like piecewise curves, we have classified them as piecewise wavy structure. 


Finally we showcase a continuous wavy structure in Fig. \ref{fig:NegativeStrengthsPattern} (b) for $\sigma = -0.0035, \mu = -0.0005, k=-1$. Since the end node plot resembles a continuous curve like structure, hence we have classified it as continuous wavy structure. 

The presence of rich spatiotemporal structures call for basin of attraction of the spatiotemporal states in the network of Chialvo neuron. \blue{This} could give us an understanding of the prevalence of chimera states, synchronized states, piecewise wavy structures, \blu{and} continuous wavy structures in the space of coupling parameters. 

\section{Conclusions and Future directions}
\label{sec:conclusions}
In this paper, we have \blu{analyzed} the dynamical effects of the inclusion of electromagnetic field on the Chialvo neuron map. We showed that the improved Chialvo map exhibits variety of rich dynamics such as many number of fixed points, multistability, antimonotonicity, reverse doubling bifurcation, presence of fingered attractors, firing patterns, revealing various global bifurcations via one and two parameter bifurcation analysis, and prevalence of chimera states in a network of Chialvo neurons. \blue{It has also been shown that both the two-dimensional Chialvo map and the three-dimensional Chialvo map in the presence of electromagnetic flux are both non-invertible smooth maps. The dynamical analysis, one-parameter, two-parameter bifurcation diagrams were also performed using software such as \sc{MATCONTM}}. This makes it a promising role model for further future research of exploring a lot of other interesting nonlinear phenomena, \blue{exploring both  local and global bifurcations analytically}.

Preliminary research \blu{has} shown that spiral waves exist in a \blu{two-dimensional} lattice network of Chialvo neurons. Many synchronization, spiral wave patterns, spiral wave chimera can be explored in two dimensional, multilayer network of Chialvo neurons according to similar works in \cite{ShMu20a,ShMu21a}. Basin of attraction of different spatiotemporal patterns can be studied which could further clarify the interesting new regimes of piecewise-wavy and continuous wavy spatiotemporal patterns. We hope this could open doors for further analytical study of bifurcations in the Chialvo neuron map. Does Chialvo map exhibit torus, torus doubling bifurcations, hyperchaos? is an open problem to investigate. \blue{It would be interesting to compute the stable and unstable manifolds of saddle fixed points near periodic orbits. This can open doors to \blu{understanding} the homoclinic theory in neuron maps.}
\section*{Acknowledgements}
The authors thank Dr. David J.W. Simpson for his teachings, supervision, techniques on dynamical systems theory. His teachings were very fundamental while exploring the interesting world of dynamical systems. His office door was always open for discussions. S.S.M acknowledges the fruitful discussions with Dr. Astero Provata regarding the feedbacks, criticisms specifically on traveling waves, spatiotemporal patterns, usage of negative coupling strengths in the network study.  The authors thank Aasifa Rounak for discussions which improved the manuscript. The authors also thank the anonymous reviewers whose comments improved the quality of the manuscript and led to new research questions. S.S.M acknowledges the School of Mathematical and Computational Sciences doctoral bursary funding during this research.  

\bibliographystyle{plain}
\bibliography{Main}

\end{document}